\documentclass[10pt]{article}
\usepackage{epsf}
\usepackage{epsfig}
\usepackage{url}

\def\RR{\rm \hbox{I\kern-.2em\hbox{R}}}
\def\NN{\rm \hbox{I\kern-.2em\hbox{N}}}
\def\ZZ{\rm {{\rm Z}\kern-.28em{\rm Z}}}
\def\CC{\rm \hbox{C\kern -.5em {\raise .32ex \hbox{$\scriptscriptstyle
|$}}\kern-.22em{\raise .6ex \hbox{$\scriptscriptstyle |$}}\kern .4em}}

\def\<{\langle}
\def\>{\rangle}
\def\t{\tilde}

\def\e{\varepsilon}

\def\nl{\newline}
\def\o{\overline}

\def\bq{{\bf q}}
\def\cT{{\cal T}}
\def\cA{{\cal A}}

\def\cB{{\cal B}}

\def\cD{{\cal D}}

\def\cC{{\cal C}}

\def\Chi{\raise .3ex
\hbox{\large $\chi$}} 
\def\lsima{\hbox{\kern -.6em\raisebox{-1ex}{$~\stackrel{\textstyle<}{\sim}~$}}\kern -.4em}
\def\lsim{\hbox{\kern -.2em\raisebox{-1ex}{$~\stackrel{\textstyle<}{\sim}~$}}\kern -.2em}
\def\({\Bigl (}
\def\){\Bigr )}
\def\({\Bigl (}
\def\){\Bigr )}

\newcommand{\be}{\begin{equation}}
\newcommand{\ee}{\end{equation}}
\newcommand{\bea}{$$ \begin{array}{lll}}
\newcommand{\eea}{\end{array} $$}
\newcommand{\bi}{\begin{itemize}}
\newcommand{\ei}{\end{itemize}}
\newcommand{\iref}[1]{(\ref{#1})}
\newtheorem{theorem}{Theorem}[section]
\newtheorem{remark}[theorem]{Remark}
\newtheorem{lemma}[theorem]{Lemma}
\newtheorem{proposition}[theorem]{Proposition}
\newtheorem{corollary}[theorem]{Corollary}
\newtheorem{definition}[theorem]{Definition}

\newtheorem{prop}[theorem]{Proposition}

\def\proof{{\noindent \bf Proof: }}

\def\ve{\varepsilon}
\def\sq{\hfill $\diamond$\\}

\def\R{\mathbb R}

\def\mE{\mathbb E}

\def\nsdf{\|\sqrt{\det(d^2f)}\|}

\def\TEq{{T_{\rm eq}}}
\def\TO{{T_0}}
\def\sep{\; ; \;}

\def\gsim{\hbox{\kern -.2em\raisebox{-1ex}{$~\stackrel{\textstyle>}{\sim}~$}}\kern -.2em}

\usepackage{amsmath}
\usepackage{amsfonts}
\usepackage{amssymb}

\usepackage{psfrag}
\usepackage{graphicx}

\DeclareMathOperator{\diam}{diam}

\def\cPi{{\text{\large$\boldsymbol\pi$}}}

\begin{document}
\title{\bf Greedy bisection generates optimally adapted triangulations} 
\author{Jean-Marie Mirebeau and Albert Cohen}
\maketitle
\date{}
\begin{abstract}
We study the properties of a simple 
greedy algorithm introduced in \cite{CDHM} for the generation of data-adapted
anisotropic triangulations. Given a function $f$, the algorithm produces nested triangulations
$\cT_N$ and corresponding piecewise polynomial
approximations $f_N$ of $f$. The refinement procedure picks 
the triangle which maximizes the local $L^p$ approximation error, and
bisects it in a direction which is chosen so to minimize this error
at the next step. We study the approximation error in the $L^p$ norm
when the algorithm is applied to $C^2$ functions
with piecewise linear approximations. We prove that 
as the algorithm progresses, the triangles tend to adopt an optimal aspect ratio which
is dictated by the local hessian of $f$. 
For convex functions, we also prove that the adaptive triangulations
satisfy the convergence bound $\|f-f_N\|_{L^p}  \leq CN^{-1}\|\sqrt{\det(d^2f)}\|_{L^\tau}$
with $\frac 1 \tau:=\frac 1 p + 1$,
which is known to be asymptotically optimal among all possible triangulations.
\end{abstract}

\section{Introduction}
\label{intro}

In finite element approximation,
a classical and important distinction is made between {\it uniform} 
and {\it adaptive} methods. In the first case
all the elements which constitute the mesh have comparable shape and size,
while these attributes are allowed to vary strongly in the second case.
An important feature of adaptive methods is the fact that 
the mesh is not fixed in advance
but rather tailored to the properties of the function $f$
to be approximated. Since the function approximating $f$ is not picked from a fixed linear space, 
adaptive finite elements can be considered as an 
instance of {\it non-linear approximation}. Other
instances include approximation by rational
functions, or by $N$-term linear combinations
of a basis or dictionary. We refer to \cite{De} for a general survey
on non-linear approximation.

In this paper, we focus our interest
on {\it piecewise linear} finite element functions defined over triangulations
of a bidimensional polygonal domain $\Omega\subset \RR^2$.  Given a triangulation
$\cT$ we denote by $V_{\cT}:=\{v\;{\rm s.t.} \; v_{|T}\in\Pi_1,\; T\in\cT\}$
the associated finite element space. The norm in which we measure the approximation error is
the $L^p$ norm for $1\leq p\leq \infty$ and we therefore do not require 
that the triangulations are conforming and that the functions of $V_\cT$ are continuous
between triangles. For a given function $f$
we define 
$$
e_N(f)_{L^p}:=\inf_{\#(\cT)\leq N} \inf_{g\in V_{\cT}} \|f-g\|_{L^p},
$$
the best approximation error of $f$ when using at most $N$ elements.
In adaptive finite element approximation, critical questions are:
\begin{enumerate}
\item
Given a function $f$ and a number $N>0$, how can we characterize the {\it optimal mesh}
for $f$ with $N$ elements corresponding to the above defined best approximation error.
\item
What quantitative estimates
are available for the best approximation error $e_N(f)_{L^p}$ ? Such estimates should involve
the derivatives of $f$ in a different way than for non-adaptive meshes.
\item
Can we build by a simple algorithmic procedure
a mesh $\cT_N$ of cardinality $N$ and a finite element function
$f_N\in V_{\cT_N}$ such that $\|f-f_N\|_{L^p}$ 
is comparable to $e_N(f)_{L^p}$ ?
\end{enumerate}

While the optimal mesh
is usually difficult to characterize exactly, it should
satisfy two intuitively desirable features:
(i) the triangulation should {\it equidistribute} the local approximation error between 
each triangle and (ii) the aspect ratio of a triangle $T$ should be 
{\it isotropic} with respect to a distorted
metric induced by the local value of the hessian $d^2f$ on $T$
(and therefore anisotropic in the sense of the euclidean metric). 
Under such prescriptions on the mesh, 
quantitative error estimates have recently
been obtained in \cite{CSX,BBLS} when $f$ is 
a $C^2$ function. These estimates
are of the form
\be
\limsup_{N\to\infty} N e_N(f)_{L^p}  \leq C\|\sqrt{|\det(d^2f)|}\|_{L^\tau},\;\; \frac 1 \tau=\frac 1 p + 1,
\label{aniser}
\ee
where $\det(d^2f)$ is the determinant of the $2\times 2$ hessian matrix.
For a convex $C^2$ function $f$ this estimate has been proved to be {\it asymptotically optimal} in \cite{CSX}, in the following sense
\be
\liminf_{N\to +\infty} Ne_N(f)_{L^p} \geq c\|\sqrt{|\det(d^2f)|}\|_{L^\tau}.
\label{lowerboundasym}
\ee
The convexity assumption can actually be replaced by a mild assumption
on the sequence of triangulations which is used for the approximation of $f$: a
sequence $(\cT_N)_{N\geq N_0}$ is said to be admissible if $\#(\cT_N)\leq N$ and
$$
\sup_{N\geq N_0}\left( N^{1/2}\max_{T\in\cT_N} \diam(T) \right)< \infty.
$$
Then it is proved in \cite{Mi}, that for any admissible sequence and any $C^2$ function $f$,
one has
\be
\liminf_{N\to +\infty} N\inf_{g\in V_{\cT_N}}\|f-g\|_{L^p} \geq c\|\sqrt{|\det(d^2f)|}\|_{L^\tau}.
\label{lowerboundadmiss}
\ee
The admissibility assumption is not a severe limitation for an upper estimate
of the error since it is also proved that for all $\e>0$, there exist an admissible sequence
such that 
\be
\limsup_{N\to +\infty} N\inf_{g\in V_{\cT_N}}\|f-g\|_{L^p} \leq C\|\sqrt{|\det(d^2f)|}\|_{L^\tau}+\e.
\label{upperboundadmiss}
\ee
We also refer to \cite{Mi} for a generalization of such upper and lower estimates
to higher order elements. 

From the computational viewpoint, a commonly used strategy
for designing an optimal mesh consists therefore in
evaluating the hessian $d^2f$ and imposing that 
each triangle of the mesh is isotropic with respect
to a metric which is properly related to its local value.
We refer in particular to \cite{BFGLS}
where this program is executed using 
Delaunay mesh generation techniques.
While these algorithms fastly produce 
anisotropic meshes which are naturally
adapted to the approximated function,
they suffer from two intrinsic limitations:
\begin{enumerate}
\item
They use the data of $d^2f$,
and therefore do not apply 
to non-smooth or noisy functions. 
\item
They are non-hierarchical: for $N>M$, the triangulation
$\cT_N$ is not a refinement of $\cT_M$. 
\end{enumerate}

In \cite{CDHM}, an alternate strategy was proposed 
for the design of adaptive hierarchical meshes, based on a 
simple {\it greedy algorithm}: starting from an initial 
triangulation $\cT_{N_0}$, the algorithm
picks the triangle $T\in \cT_N$ with the largest local
$L^p$ error. This triangle is then bisected
from the mid-point of one of its edges to the opposite vertex. The choice of the edge
among the three options is 
the one that minimizes the new approximation
error after bisection.
The algorithm can be applied to any $L^p$ function,
smooth or not, in the context of 
piecewise polynomial approximation of any given order.
In the case of piecewise linear
approximation, numerical experiments 
in \cite{CDHM} indicate that this elementary strategy
generates triangles with an optimal aspect ratio
and approximations $f_N\in V_{\cT_N}$ such that $\|f-f_N\|_{L^p}$
satisfies the same estimate as $e_N(f)_{L^p}$ in \iref{aniser}.

The goal of this paper is to support these experimental observations
by a rigorous analysis. Our paper is organized as follows:

In \S 2, we introduce 
notations which are used throughout the paper
and collect some available approximation
theory results for piecewise linear finite elements,
making the distinction between (i) uniform,
(ii) adaptive isotropic and (iii) adaptive
anisotropic triangulations. In the last case, which
is in the scope of this paper, we introduce a measure
of non-degeneracy of a triangle $T$ with respect
to a quadratic form. We show that the optimal
error estimate \iref{aniser} is met when each triangle
is non-degenerate in the sense of the above
measure with respect to the quadratic form
given by the local hessian $d^2f$.
We end by briefly recalling the greedy algorithm
which was introduced in \cite{CDHM}. 

In \S 3, we study the behavior 
of the refinement procedure when applied to a
quadratic function $q$ such that its associated quadratic
form $\bq$ is of positive or negative sign.
A key observation is that the edge which is bisected
is the longest with respect to the metric induced by $\bq$.
This allows us to prove that the triangles generated 
by the refinement procedure adopt 
an optimal aspect ratio in the sense
of the non-degeneracy measure introduced in \S 2.

In \S 4, we study the behavior
of the algorithm when applied to a general $C^2$
function $f$ which is assumed to be
strictly convex (or strictly concave). We first establish a perturbation
result, which shows that when $f$ is locally
close to a quadratic function $q$ the algorithm
behaves in a similar manner as when applied
to $q$. We then prove that the diameters of
the triangles produced by the algorithm tend
to zero so that the perturbation result can be 
applied. This allows us to show that the
optimal convergence estimate 
\be
\limsup_{N \to \infty} N \|f-f_N\|_{L^p} \leq C \|\sqrt{|\det(d^2f)|}\|_{L^\tau}
\label{opticonvest}
\ee
is met by the sequence of approximations
$f_N\in V_{\cT_N}$ generated by the algorithm.

The extension of this result to an arbitrary $C^2$ function
$f$ remains an open problem. It is possible to
proceed to an analysis similar to \S 3 in the case
where the quadratic form $\bq$ is of mixed sign, also proving that
the triangles adopt an optimal aspect ratio as they get refined.
We describe this analysis in \S 9.1 of \cite{thesisJM}. However,
it seems difficult to extend the perturbation
analysis of \S 4 to this new setting. In particular
the diameters of the triangles are no more ensured to
tend to zero, and one can even exhibit examples of non-convex $C^2$
functions $f$ for which the approximation $f_N$ {\it fails} to converge towards $f$
due to this phenomenon. Such examples are
discussed in \cite{CDHM} which also proposes
a modification of the algorithm for which convergence
is always ensured. However, we do not know if
the optimal convergence estimate \iref{opticonvest} holds for
any $f\in\cC^2$ with this modified algorithm, although this seems
plausible from the numerical experiments.

\section{Adaptive finite element approximation}

\subsection{Notations}

We shall make use of a linear approximation operator
$\cA_T$ that maps continuous functions defined on $T$ onto $\Pi_1$. 
For an arbitrary but fixed  $1\leq p\leq \infty$, we define the
local $L^p$ approximation error 
$$
e_T(f)_p:=\|f-\cA_Tf\|_{L^p(T)}.
$$
The critical assumptions in our analysis for the operator $\cA_T$ 
will be the following:
\begin{enumerate}
\item
$\cA_T$ is continuous in the $L^\infty$ norm. 
\item
$\cA_T$ commutes with affine changes of variables: $\cA_T(f) \circ \phi = \cA_{\phi^{-1}(T)} (f\circ \phi)$
for all affine $\phi$.
\item
$\cA_T$ reproduces $\Pi_1$:
$\cA_T(\pi)=\pi$,
for any $\pi\in\Pi_1$.
\end{enumerate}
Note that the commutation assumption implies that for any function $f$
and any affine transformation $\phi:x\mapsto x_0+Lx$ we have
\be
e_{\phi(T)}(f)_p=|\det(L)|^{1/p} e_T(f\circ \phi)_p,
\label{commuterror}
\ee
Two particularly simple admissible choices of approximation operators are the following:
\begin{itemize}
\item
$\cA_T=P_T$, the $L^2(T)$-orthogonal projection operator: $\int_T (f-P_T f)\pi=0$ for all $\pi\in \Pi_1$.
\item
$\cA_T=I_T$, the local interpolation operator: 
$I_Tf(v_i)=f(v_i)$ with $\{v_0,v_1,v_2\}$ the vertices of $T$.
\end{itemize}
All our results are simultaneously valid when $\cA_T$ is either $P_T$
or $I_T$, or any linear operator that fulfills the three above assumptions.

Given a function $f$ and a triangulation 
$\cT_N$ with $N=\#(\cT_N)$, we can associate a finite element approximation $f_N$
defined on each $T\in\cT_N$ by
$f_N(x)=\cA_T f(x)$. The global approximation error is  given by
$$
\|f-f_N\|_{L^p}=\(\sum_{T\in\cT_N}e_T(f)_p^p\)^{\frac 1 p},
$$
with the usual modification when $p=\infty$. 

\begin{remark}
The operator $\cB_T$ of best $L^p(T)$ approximation
which is defined by
$$
\|f-\cB_Tf\|_{L^p(T)}=\min_{\pi\in\Pi_m}\|f-\pi\|_{L^p(T)},
$$ 
does not fall in the above category of operators, since
it is non-linear (and not easy to compute) when $p\neq 2$. However, it is clear
that any estimate on $\|f-f_N\|_{L^p}$ with $f_N$ defined as $\cA_T f$ on each $T$
implies a similar estimate when $f_N$ is defined as $\cB_T f$ on each $T$.
\end{remark}

Here and throughout the paper, when
$$
q(x,y)= a_{2,0} x^2+ 2 a_{1,1} xy+ a_{0,2} y^2+ a_{1,0} x+ a_{0,1} y + a_{0,0}
$$
we denote by $\bq$ the associated quadratic form : if $u=(x,y)$
$$
\bq(u)= a_{2,0} x^2+ 2 a_{1,1} xy+ a_{0,2} y^2.
$$
Note that $\bq(u)=\<Qu,u\>$ 
where 
$
Q=\left(\begin{array}{cc}
a_{2,0} & a_{1,1}\\
a_{1,1} & a_{0,2}
\end{array} \right).
$
We define
$$
\det (\bq) := \det (Q).
$$
If $\bq$ is a positive or negative quadratic form, we define the
{\it $\bq$-metric}
\be
|v|_\bq:=\sqrt{|\bq(v)|}
\label{qmetric}
\ee
which coincides with the euclidean norm when $\bq(v)=x^2+y^2$ for $v=(x,y)$.
If $\bq$ is a quadratic form of mixed sign, we define the associated
positive form $|\bq|$ which corresponds to the symmetric matrix $|Q|$
that has same eigenvectors as $Q$ with eigenvalues $(|\lambda|,|\mu|)$
if $(\lambda,\mu)$ are the eigenvalues of $Q$. Note that generally $|\bq|(u)\neq |\bq(u)|$
and that one always has $|\bq(u)|\leq |\bq|(u)$.

\begin{remark}
\label{rmCanonical}
If $\det Q>0$, then there exists a $2\times 2$ matrix $L$ and $\ve\in \{+1, -1\}$ such that 
$$
L^t Q L = \ve 
\left(
\begin{array}{cc}
1 & 0\\
0 & 1
\end{array}
\right).
$$
The linear change of coordinates $\phi(u) := Lu$, where $u = (x,y) \in \R^2$, therefore satisfies $\bq \circ \phi (u) = \ve (x^2+y^2)$. On the other hand, if $\det Q<0$ then there exists a $2 \times 2$ matrix $L$ such that 
$$
L^t Q L = 
\left(
\begin{array}{cc}
1 & 0\\
0 & -1
\end{array}
\right).
$$
Defining again $\phi(u) := Lu$ we obtain in this case $\bq \circ \phi (u) = x^2 - y^2$.
\end{remark}
\subsection{From uniform to adaptive isotropic triangulations}

A standard estimate in finite element approximation states that 
if $f\in W^{2,p}(\Omega)$ then
$$
\inf_{g\in V_h} \|f-g\|_{L^p} \leq C h^{2}\|d^2f\|_{L^{p}},
$$
where $V_h$ is the piecewise linear finite element space
associated with a triangulation $\cT_h$ of mesh size
$h:=\max_{T\in\cT_h} \diam(T)$. If we restrict our attention
to {\it uniform} triangulations, we have
$$
N:=\#(\cT_h) \sim  h^{-2}.
$$
Therefore, denoting by $e^{\rm unif}_N(f)_{L^p}$ the $L^p$ approximation error
by a uniform triangulation of cardinality $N$, we can re-express the above
estimate as
\be
\label{uniferror}
e^{\rm unif}_N(f)_{L^p} \leq CN^{-1}\|d^2f\|_{L^p}.
\ee
This estimate can be significantly improved when using
adaptive partitions. We give here some heuristic
arguments, which are based on the assumption that on each triangle $T$
the relative variation of $d^2f$ is small so that it can be considered
as constant over $T$ (which means that 
$f$ is replaced by a quadratic function on each $T$),
and we also indicate the available results which 
are proved more rigorously.

First consider {\it isotropic} triangulations,
i.e. such that all triangles satisfy a uniform estimate 
\be
\rho_T=\frac {h_T}{r_T}\leq A,
\label{regulartri}
\ee
where $h_T:=\diam(T)$ denotes the size of the longest edge of $T$, and $r_T$ is the radius of the largest disc contained in $T$. 
In such a case we start from the local approximation estimate on any $T$
$$
e_T(f)_p \leq Ch_T^2 \|d^2f\|_{L^p(T)},
$$
and notice that 
$$
h_T^2 \|d^2f\|_{L^p(T)}\sim |T|\, \|d^2f\|_{L^p(T)} = \|d^2f\|_{L^\tau(T)},
$$
with $\frac 1 \tau:=\frac 1 p +1$ and $|T|$ the area of $T$, where we have used
the isotropy assumption \iref{regulartri} in the equivalence and the fact that $d^2f$
is constant over $T$ in the equality. It follows that
$$
e_T(f)_p \leq C \|d^2f\|_{L^\tau(T)}, \;\; \frac 1 \tau:=\frac 1 p +1.
$$
Assume now that we can construct adaptive isotropic triangulations $\cT_N$ 
with $N:=\#(\cT_N)$ which {\it equidistributes} the local error in the
sense that for some prescribed $\e>0$
\be
c\e \leq e_T(f)_p  \leq \e,
\label{equil}
\ee
with $c>0$ a fixed constant independent of $T$ and $N$. Then
defining $f_N$ as $\cA_T(f)$ on each $T\in\cT_N$, we have on the one hand
$$
\|f-f_N\|_{L^p} \leq N^{1/p}\e,
$$
and on the other hand, with $\frac 1 \tau:=\frac 1 p +1$,
$$
N(c\e)^\tau \leq \sum_{T\in\cT_N}\|f-f_N\|_{L^p(T)}^\tau 
\leq C^\tau\sum_{T\in\cT_N} \|d^2f\|_{L^\tau(T)}^\tau\leq C^\tau \|d^2f\|_{L^\tau}^\tau.
$$
Combining both, one obtains for $e^{\rm iso}_N(f)_{L^p}:=\|f-f_N\|_{L^p}$
the estimate
\be
\label{isoerror}
e^{\rm iso}_N(f)_{L^p} \leq  CN^{-1}\|d^2f\|_{L^\tau}.
\ee
This estimate improves upon \iref{uniferror} since the rate $N^{-1}$ is now
obtained with the weaker smoothness condition $d^2f\in L^\tau$
and since, even for smooth $f$, the quantity $\|d^2f\|_{L^\tau}$ might
be significantly smaller than $\|d^2f\|_{L^p}$.
This type of result is classical in non-linear approximation and also occurs
when we consider best $N$-term approximation in
a wavelet basis.

The principle of error equidistribution suggests a
simple {\it greedy algorithm} to build an adaptive isotropic triangulation 
for a given $f$, similar to our algorithm but
where the bisection of the triangle $T$
that maximizes
the local error $e_T(f)_p$ is systematically done from
its {\it most recently created vertex} in order to preserve the estimate
\iref{regulartri}. Such an algorithm cannot exactly equilibrate the error
in the sense of \iref{equil} and therefore does not lead to the 
same the optimal estimate as in \iref{isoerror}.
However, it was proved in \cite{BDDP}  that it satisfies
$$
\|f-f_N\|_{L^p} \leq C |f|_{B^2_{\tau,\tau}}N^{-1},
$$
for all $\tau$ such that $\frac 1 \tau < \frac 1 p +1$, provided that the
local approximation operator $\cA_T$ is bounded in the $L^p$ norm. 
Here $B^2_{\tau,\tau}$ denotes the usual Besov space which is a 
natural substitute for $W^{2,\tau}$ when $\tau<1$. Therefore
this estimate is not far from \iref{isoerror}.

\subsection{Anisotropic triangulations: the optimal aspect ratio}

We now turn to anisotropic adaptive triangulations, and start by
discussing the optimal shape of a triangle $T$ for a given function $f$
at a given point. For this purpose,
we again replace $f$ by a quadratic function
assuming that $d^2f$ is constant over $T$. For such
a $q\in \Pi_2$ and its associated quadratic form $\bq$, we 
first derive an equivalent quantity for the local approximation error.
Here and as well as in \S 3 and \S 4, we consider a triangle
$T$ and we denote by $(a,b,c)$
its edge vectors oriented in clockwise 
or anti-clockwise direction so that
$$
a+b+c=0.
$$
\begin{prop}
\label{propequiverrornondeg}
The local $L^p$-approximation error satisfies
$$
e_T(q)_p=e_T(\bq)_p\sim |T|^{\frac 1 p} \max\{|\bq(a)|,|\bq(b)|,|\bq(c)|\},
$$
where the constant in the equivalence is independent of $q$, $T$ and $p$.  
\end{prop}

\proof
The first equality is trivial since $q$ and $\bq$ differ by an affine function.
Let $\TEq$ be an equilateral triangle of area $|\TEq|=1$, and edges $a,b,c$. 
Let $E$ be the $3$-dimensional vector space of all quadratic forms.
Then the following quantities are norms on $E$, and thus equivalent: 
\be
\label{equivTEq}
e_\TEq(\bq)_p \sim \max\{|\bq(a)|,|\bq(b)|,|\bq(c)|\}.
\ee
Note that the constants in this equivalence are independent of $p$ since all $L^p(T)$ norms
are uniformly equivalent on $E$. 

If $T$ is an arbitrary triangle, there exists an affine transform $\phi: x\mapsto x_0+Lx$ such that $T =\phi(\TEq)$.
For any quadratic function $q$, we thus obtain from \iref{commuterror}
$$
e_T(\bq)=e_T(q) = e_{\phi(\TEq)} (q) = |\det L|^{\frac 1 p} e_\TEq (q\circ \phi)= |\det L|^{\frac 1 p} e_\TEq (\bq\circ L)
$$
since $\bq \circ L$ is the homogeneous part of $q\circ \phi$.
By \iref{equivTEq}, we thus have
$$
e_T(\bq)
 \sim |\det L|^{\frac 1 p} \max\{|\bq(La)|,\ |\bq(Lb)|,\ |\bq(Lc)|\},
$$
where $\{a,b,c\}$ are again the edge vectors of $\TEq$.
Remarking that $|T| = |\det L|$ and that $\{La,Lb,Lc\}$ are the edge vectors of $T$, 
this concludes the proof of this proposition.
\sq
\nl
In order to describe the optimal shape of a triangle $T$ for the quadratic function
$q$, we fix the area of $|T|$ and try to minimize the error  $e_T(q)_p$
or equivalently $\max\{|\bq(a)|,|\bq(b)|,|\bq(c)|\}$. The solution to this 
problem can be found by introducing for any $\bq$ such that 
$\det(\bq)\neq 0$ the following measure of {\it non-degeneracy} for $T$:
\be
\label{rhodef}
\rho_\bq(T):=\frac {\max\{|\bq(a)|,|\bq(b)|,|\bq(c)|\}}{|T|\sqrt{|\det(\bq)|}}.
\ee
Let $\phi$ be a linear change of variables, $\bq$ a quadratic form and $T$ a triangle of edges $a,b,c$. Then $\det(\bq\circ \phi) = (\det \phi)^2 \det (\bq)$, the edges of $\phi(T)$ are $\phi(a), \phi(b), \phi(c)$ and $|\phi(T)| = |\det \phi| |T|$.
Hence we obtain
\be
\label{rhoinv}
\rho_{\bq\circ \phi}(T)
= \frac {\max\{|\bq\circ \phi (a)|,|\bq \circ \phi (b)|,|\bq \circ \phi (c)|\}}{|T|\sqrt{|\det(\bq \circ \phi)|}} 
= \frac {\max\{|\bq(\phi(a))|,|\bq(\phi(b))|,|\bq(\phi(c))|\}}{|\det\phi| |T|\sqrt{|\det(\bq)|}} 
= \rho_\bq(\phi(T)).
\ee


\noindent
The last equation, combined with Remark \ref{rmCanonical}, allows to reduce the study of $\rho_\bq(T)$ to two 
elementary cases by change of variable:
\begin{enumerate}
\item
The case where $\det(\bq)>0$ is reduced to
$\bq(x,y)=x^2+y^2$. 
Recall that for any triangle $T$ with edges $a,b,c$ we define $h_T := \diam(T) = \max\{|a|,|b|,|c|\}$,
with $|\cdot|$ the euclidean norm.
In this case we therefore have $\rho_\bq(T)=\frac {h_T^2}{|T|}$,
which corresponds to a standard measure of shape regularity in the sense that its boundedness
is equivalent to a property such as \iref{regulartri}. 
This quantity is minimized when the triangle $T$ is equilateral, with
minimal value $\frac 4{\sqrt 3}$ (in fact it was also proved in \cite{Chen1}
that the minimum of the interpolation error
$\|\bq-I_T\bq\|_{L^p(T)}$ among all triangles of area $|T|=1$ is attained 
when $T$ is equilateral). For a general quadratic form $\bq$ of positive
sign, we obtain by change of variable that the minimal value $\frac  4{\sqrt 3}$
is obtained for triangles which are equilateral with respect to the 
metric $|\cdot|_{\bq}$. More generally triangles with a good aspect ratio,
i.e. a small value of $\rho_{\bq}(T)$,
are those which are {\it isotropic with respect to this metric}. Of course, a similar
conclusion holds for a quadratic form of negative sign. 
\item
The case where $\det(\bq)<0$ is reduced to
$\bq(x,y)=x^2-y^2$. In this case, the analysis presented in \cite{Cao}
shows that the quantity $\rho_\bq(T)$ is minimized when
$T$ is a half of a square with sides parallel to the $x$ and $y$ axes,
with minimal value $2$. But using \iref{rhoinv} we also notice that 
$\rho_\bq(T) = \rho_\bq (L(T))$ for any linear transformation $L$
such that $\bq = \bq\circ L$. This holds if $L$ has eigenvalues $(\lambda,\frac 1 \lambda)$, 
where $\lambda \neq 0$, and eigenvectors $(1,1)$ and $(-1,1)$. 
Therefore, all images of the half square by such transformations
$L$ are also optimal triangles. Note that such triangles can be highly anisotropic.
For a general quadratic form $\bq$ of mixed sign, we notice that $\rho_\bq(T)\leq \rho_{|\bq|}(T)$, 
and therefore triangles which are equilateral
with respect to the metric $|\cdot|_{|\bq|}$ have 
a good aspect ratio, i.e. a small value of $\rho_\bq(T)$. 
In addition, by similar arguments, we find that all images
of such triangles by linear transforms $L$ with eigenvalues
$(\lambda,\frac 1 \lambda)$ and eigenvectors $(u,v)$ such that $\bq(u)=\bq(v)=0$
also have a good aspect ratio, since 
$\bq = \bq\circ L$ for such transforms.
\end{enumerate}
We leave aside the special case where $\det (\bq)=0$. In such a case, 
the triangles minimizing the error for a given area degenerate
in the sense that they should be infinitely long and thin,
aligned with the direction of the null eigenvalue of $\bq$.

Summing up, we find that triangles with a good aspect ratio
are characterized by the fact that $\rho_\bq(T)$ is small.
In addition, from Proposition \ref{propequiverrornondeg} and the definition
of $\rho_\bq(T)$, we have
\be
e_T(q)_p\sim |T|^{1+\frac 1 p}\sqrt{|\det(\bq)|}\rho_\bq(T) = \|\sqrt{|\det(\bq)|}\|_{L^\tau(T)}
\rho_\bq(T),\;\;  \frac 1 \tau:=\frac 1 p+1.
\label{errorrhoq}
\ee
We now return to a function $f$ such that $d^2f$ is assumed to be constant
on every $T\in\cT_N$. Assuming that all triangles have a good
aspect ratio in the sense that
$$
\rho_\bq(T) \leq C
$$
for some fixed constant $C$ and with $\bq$ the value of $d^2f$ over $T$, we find
up to a change in $C$ that 
\be
e_T(f)_p\leq C\|\sqrt{|\det(d^2f)|}\|_{L^\tau(T)}
\ee
By a similar reasoning as with isotropic triangulations, we now obtain
that if the triangulation equidistributes the error in the sense of \iref{equil}
\be
\label{aniserror2}
\|f-f_N\|_{L^p} \leq  CN^{-1}\|\sqrt{|\det(d^2f)|}\|_{L^\tau},
\ee
and therefore \iref{aniser} holds. This estimate improves
upon \iref{isoerror} since the quantity $\|\sqrt{|\det(d^2f)|}\|_{L^\tau}$
might be significantly smaller than $\|d^2f\|_{L^\tau}$, in particular
when $f$ has some anisotropic features, such as sharp gradients
along curved edges.

The above derivation of \iref{aniser} is heuristic and non-rigorous.
Clearly, this estimate cannot be valid as such since 
$\det(d^2f)$ may vanish while the approximation error does not
(consider for instance $f$ depending only on a single variable). 
More rigorous versions were derived in 
\cite{CSX} and \cite{BBLS}. 
In these results $|d^2f|$ is typically replaced by a majorant $|d^2f|+\e I$,
avoiding that its determinant vanishes.
The estimate \iref{aniser} can then be rigorously proved but 
holds for $N\geq N(\e,f)$ large enough. This limitation
is unavoidable and reflects the
fact that enough resolution is needed so that the hessian
can be viewed as locally constant over each optimized triangle.
Another formulation,
which is rigorously proved in \cite{Mi}, reads as follows.

\begin{prop}
There exists an absolute constant $C>0$ such that for any polygonal domain $\Omega$ and any function $f\in C^2(\overline \Omega)$, one has
$$
\limsup_{N\to +\infty} Ne_N(f)_{L^p} \leq C\|\sqrt{|\det(d^2f)|}\|_{L^\tau}.
$$
 \end{prop}
\subsection{The greedy algorithm}

Given a target function $f$, our algorithm iteratively builds triangulations
$\cT_N$ with $N=\#(\cT_N)$ and finite element approximations $f_N$.
The starting point is
a coarse triangulation $\cT_{N_0}$. Given $\cT_N$, the algorithm
selects the triangle $T$ which maximizes the local error $e_T(f)_p$
among all triangles of $\cT_N$, and bisects it from the mid-point of one of its edges towards the opposite vertex.
This gives the new triangulation $\cT_{N+1}$.

The critical part of the algorithm lies in the choice of the edge
$e\in\{a,b,c\}$ from which $T$ is bisected. Denoting by
$T_e^1$ and $T_e^2$
the two resulting triangles, we choose $e$
as the minimizer of  a {\it decision function} $d_T(e,f)$, which role is
to drive the generated triangles towards an optimal aspect ratio.
While the most natural choice for $d_T(e,f)$ corresponds to the split
that minimizes the error after bisection, namely
$$
d_T(e,f)=e_{T_e^1}(f)_p^p+e_{T_e^2}(f)_p^p,
\label{optisplit}
$$
we shall instead focus our attention on a decision function which 
is defined as the $L^1$ norm of the interpolation error
\be
d_T(e,f)=\|f-I_{T_e^1}f\|_{L^1(T_e^1)}+\|f-I_{T_e^2}f\|_{L^1(T_e^2)}.
\label{optil1}
\ee
For this decision, the analysis of the 
algorithm is made simpler, due to the fact that
we can derive explicit expressions of $\|f-I_T f\|_{L^1(T)}$ when $f=q$ is a quadratic polynomial
with a positive homogeneous part $\bq$. We
prove in \S 3 that this choice leads to triangles with an optimal aspect ratio
in the sense of a small $\rho_\bq(T)$. This leads us in \S 4 to a proof that
the algorithm satisfies the optimal convergence estimate \iref{aniserror2}
in the case where $f$ is $C^2$ and strictly convex.

\begin{remark}
It should be well understood that while the decision function is based on the $L^1$ norm, the
selection of the triangle to be bisected is done by maximizing $e_T(f)_p$. The algorithm remains
therefore governed by the $L^p$ norm 
in which we wish to minimize the error $\|f-f_N\|_p$ for a given number of triangles. Intuitively, this
means that the $L^p$-norm influences the size of the triangles which have to equidistribute
the error, but not their optimal shape.
\end{remark}

\begin{remark}
 It was pointed out to us that
the $L^1$ norm of the interpolation error to a suitable convex function is 
also used to improve the mesh in the context of moving grid techniques, see \cite{Chen2}.
\end{remark}

\noindent
We define a variant of the decision function as follows 
$$
D_T(e,f) := \|f-I_Tf\|_{L^1(T)} - d_T(e,f).
$$
Note that $D_T(e,f)$ is the reduction of the $L^1$ interpolation error 
resulting from the bisection of the edge $e$, and that the selected edge that minimizes $d_T(\cdot,f)$ is also the one that maximizes $D(\cdot,f)$. 
The function $D_T$ has a simple expression in the case where 
$f$ is a convex function.
\begin{lemma}
Let $T$ be a triangle and let $f$ be a  convex function on $T$. Let $e$ be an edge of $T$ with endpoints $z_0$ and $z_1$. 
Then 
\be
\label{DDiff}
D_T(e,f) = \frac {|T|} 3  \left(\frac{f(z_0)+f(z_1)} 2-f\left(\frac {z_0+z_1} 2 \right)\right).
\ee
If in addition $f$ has $C^2$ smoothness, we also have
\be
\label{Dint01}
D_T(e,f) = \frac {|T|} 6 \int_0^1 \<d^2 f(z_t) e, e\> \min\{t,1-t\} dt, \text{ where } z_t := (1-t)z_0 + t z_1.
\ee

\end{lemma}

\proof
Since $f$ is convex, we have $I_T f\geq f$ on $T$, hence 
$$
 \|f-I_Tf\|_{L^1(T)}  = \int_T (I_T f-f).
$$
Similarly $I_{T_e^1} f\geq f$ on $T_e^1$ and $I_{T_e^2} f\geq f$ on $T_e^2$, hence 
$$
D_T(e,f) = \int_T I_T f - \int_{T_e^1} I_{T_e^1} f - \int_{T_e^2} I_{T_e^2} f.
$$
Let $z_2$ be the vertex of $T$ opposite the edge $e$. 
Since the function $f$ is convex, it follows the previous expression that 
$D_T(e,f)$ is the volume of the tetrahedron of vertices 
$$
\left(\frac{z_{0,x}+ z_{1,x}} 2 ,\frac{ z_{0,y}+ z_{1,y}} 2 , f\left(\frac{z_0+z_1} 2\right)\right) \text{ and } (z_{i,x},z_{i,y}, f(z_i)) \text{ for } i=0,1,2.
$$
where $(z_{i,x},z_{i,y})$ are the coordinates of $z_i$.
Let $u = z_0-z_2$ and $v=z_1-z_2$. 
We thus have $D_T(e,f) = \frac 1 6 |\det(M)|$ where 
$$
M:= \left(
\begin{array}{ccc}
u_x & v_x & \frac {u_x+v_x} 2\\
u_y & v_y & \frac {u_y+v_y} 2\\
f(z_0) - f(z_2) & f(z_1)-f(z_2) & f\left(\frac {z_0+z_1} 2 \right) - f(z_2)
\end{array} 
\right).
$$
Subtracting the half of the first two columns to the third one we find that $M$ has the same
determinant as
$$
\t M:=\left(
\begin{array}{ccc}
u_x & v_x & 0\\
u_y & v_y & 0\\
f(z_0) - f(z_2) & f(z_1)-f(z_2) & f\left(\frac {z_0+z_1} 2 \right) - \frac{f(z_0)+f(z_1)} 2
\end{array} 
\right).
$$
Recalling that $2|T| = |\det(u,v)|$ we therefore obtain 
\iref{DDiff}. In order to
establish \iref{Dint01},
we observe that we have in the distribution sense $\partial_t^2 (\min\{t,1-t\}^+) = \delta_0-2 \delta_{1/2}+ \delta_1$, 
where $\delta_t$ is the one-dimensional Dirac function at a point $t$. Hence for any univariate function $h\in C^2([0,1])$, we have 
$$
\int_0^1 h''(t) \min\{t,1-t\} dt = h(0)-2h(1/2)+ h(1).
$$
Combining this result with \iref{DDiff} we obtain \iref{Dint01}.
\sq

\section{Positive quadratic functions}

In this section, we study the algorithm 
when applied to a quadratic polynomial $q$ such that $\det(\bq)>0$. 
We shall assume without loss of generality that
$\bq$ is positive definite, since all our results
extend in a trivial manner to the negative definite case.

Our first observation is that 
the refinement procedure based on the 
decision function \iref{optil1} 
always selects for bisection the {\it longest edge} in the sense of the $\bq$-metric $|\cdot|_{\bq}$
defined by \iref{qmetric}.

\begin{lemma}
An edge $e$ of $T$ maximizes $D_T(e,q)$ among all edges of $T$
if and only if it maximizes $|e|_{\bq}$ among all edges of $T$.
\end{lemma}

\proof
The hessian $d^2q$ is constant and for all $e\in \R^2$ one has 
$$
\<d^2 q e,e\> = 2\bq(e).
$$
If $e$ is an edge of a triangle $T$, and if $q$ is a convex quadratic function, equation \iref{Dint01} 
therefore gives
\be
\label{eqDq}
D_T(e,q) = \frac {|T|} {3} \bq(e) \int_0^1\min\{t,1-t\}dt= \frac {|T|} {12} |e|_\bq^2.
\ee
This concludes the proof. 
\sq

It follows from this lemma that the longest edge of $T$ in the sense of the $\bq$-metric
is selected for bisection by the decision function. In the remainder of this
section, we use this fact to prove that the refinement procedure
produces triangles which tend to adopt an
optimal aspect ratio in the sense that $\rho_\bq(T)$ becomes small
in an average sense.

For this purpose, it is convenient to introduce a close variant to 
$\rho_{\bq}(T)$: if $T$ is a triangle with edges $a,b,c$, such that $|a|_\bq\geq |b|_\bq \geq |c|_\bq$, 
we define 
\be
\sigma_\bq(T):= \frac{\bq(b)+\bq(c)}{4 |T| \sqrt{\det \bq}}= \frac{|b|_\bq^2+|c|_\bq^2}{4 |T| \sqrt{\det \bq}}.
\label{sigmaq}
\ee
Using the inequalities $|b|^2_\bq+|c|^2_\bq\leq 2|a|^2_\bq$ and 
$|a|^2_\bq \leq 2(|b|^2_\bq+|c|^2_\bq)$, we obtain the equivalence
\be
 \frac {\rho_\bq(T)} 8 \leq \sigma_\bq(T)\leq \frac {\rho_\bq(T)} 2.
 \label{equivnondeg}
\ee
Similar to $\rho_\bq$, this quantity is invariant under a linear coordinate changes $\phi$,
in the sense that
$$
\sigma_{\bq\circ \phi}(T) = \sigma_\bq (\phi(T)),
$$
From \iref{errorrhoq} and \iref{equivnondeg}
we can relate $\sigma_\bq$ to the local
approximation error.
\begin{prop}
\label{propequiverrornondeg1}
There exists a constant $C_0$, which depends only on the choice of $\cA_T$, such that for any triangle $T$, quadratic function $q$ and exponent $1 \leq p \leq \infty$, 
the local $L^p$-approximation error satisfies
\be
\label{eqEQSigma}
C_0^{-1} e_T(q)_p \leq \sigma_\bq(T) \|\sqrt{\det \bq}\|_{L^\tau(T)} \leq C_0 e_T(q)_p.
\ee
where $\frac 1 \tau:=\frac 1 p+1$.
\end{prop}
Our next result shows that $\sigma_\bq(T)$ is always reduced by 
the refinement procedure.

\begin{prop} \label{PropSigmaDec}
If $T$ is a triangle with children $T_1$ and $T_2$ obtained by the
refinement procedure for the quadratic function $q$, then
$$
\max\{\sigma_\bq(T_1), \sigma_\bq(T_2)\}\leq \sigma_\bq(T).
$$
\end{prop}

\proof
Assuming that $|a|_\bq\geq |b|_\bq\geq |c|_\bq$, we know that the edge $a$ is cut
and that the children have area $|T|/2$ and edges $a/2,b,(c-b)/2$ and $a/2, (b-c)/2,c$ 
(recall that $a+b+c=0$). We then have
\begin{eqnarray}
	2 |T| \sqrt{\det \bq}\ \sigma_\bq(T_i) & \leq & \bq\left(\frac a 2\right)+\bq\left(\frac{b-c}{2}\right)\\
	& = & \bq\left(\frac{b+c}{2}\right)+\bq\left(\frac{b-c}{2}\right)\\
	& = & \frac{\bq(b)+\bq(c)}{2}\\
	& = & 2 |T|  \sqrt{\det \bq}\ \sigma_\bq(T).
	\end{eqnarray}
\sq

\begin{remark}
When $\bq$ is the euclidean metric, the triangle that minimizes $\sigma_\bq$
is the half square. This is consistent with the above result since it
is the only triangle which is similar (i.e. identical 
up to a translation, a rotation and a dilation) to both of its children after one step of longest edge 
bisection. 
\end{remark}

\begin{remark}
A result of similar nature was already proved in \cite{Ri} : 
longest edge bisection
has the effect that the minimal angle in any triangle 
after an arbitrary number of refinements is at most twice
the minimal angle of the initial triangle.
\end{remark}

Our next objective is to show that as we iterate the refinement process,
the value of $\sigma_\bq(T)$ 
becomes bounded independently of $q$
for almost all generated triangles. For this purpose we introduce the following notation:
if $T$ is a triangle with edges such that $|a|_\bq \geq |b|_\bq \geq |c|_\bq$,
we denote by $\psi_\bq(T)$ the subtriangle of $T$ obtained after bisection of $a$
which contains the smallest edge $c$.
We first establish inequalities between the measures $\sigma_\bq$ and $\rho_\bq$
applied to $T$ and $\psi_\bq(T)$.

\begin{prop}
Let $T$ be a triangle, then
\begin{eqnarray}
\label{sigma58rho}
\sigma_\bq(\psi_\bq(T)) & \leq & \frac 5 8 \rho_\bq(T)\\
\label{rhodec}
\rho_\bq(\psi_\bq(T))   & \leq & \frac{\rho_\bq(T)} 2 \left(1+\frac{16}{\rho_\bq^2(T)}\right)
\end{eqnarray}
\end{prop}

\proof
We first prove \iref{sigma58rho}. Obviously,
$\psi_\bq(T)$ contains one edge $s\in \{a,b,c\}$ from $T$, 
and one half edge $t\in \{\frac a 2,\frac b 2,\frac c 2\}$ from $T$. Therefore
$$
\sigma_\bq(\psi_\bq(T)) \leq \frac{|s|_\bq^2+|t|_\bq^2}{4|\psi_\bq(T)|\sqrt{\det \bq}} \leq \frac{|a|_\bq^2+|\frac a 2|_\bq^2}{2|T|\sqrt{\det \bq}} = \frac 5 8 \rho_\bq(T).
$$
For the proof of \iref{rhodec}, we restrict our attention
to the case $\bq= x^2+y^2$, without loss of generality
thanks to the invariance formula \iref{rhoinv}.
Let $T$ be a triangle with edges $|a|\geq |b|\geq |c|$.
If $h$ is the width of $T$ in the direction perpendicular to $a$, then 
$$
h = \frac{2|T|}{|a|} = \frac{2 |a|}{\rho_\bq(T)}.
$$
The sub-triangle $\psi_\bq(T)$ of $T$ has edges $\frac a 2, c , d$ where $d=\frac{b-c} 2$, and the angles at the ends of $\frac a 2$ are acute. Indeed
$$
\<c,a/2\> = \frac 1 4 \left(|b|^2 - |a|^2-|c|^2\right)\leq 0 \text{ and } \<d,a/ 2\> = \frac 1 4 \left(|c|^2-|b|^2\right)\leq 0. 
$$
By Pythagora's theorem we thus find
$$
\max\{\left|\frac a 2\right|^2, |c|^2,|d|^2\} \leq \left|\frac a 2\right|^2 +h^2 = \frac{|a|^2} 4 \left(1+\frac{16}{\rho_\bq^2(T)}\right).
$$ 
Dividing by the respective areas of $T$ and $\psi_\bq(T)$, we obtain the announced result.
\sq

Our next result shows that a significant reduction of $\sigma_\bq$ occurs at least for 
one of the triangles obtained by three successive refinements, unless it has reached a
small value of $\sigma_\bq$. We use the notation $\psi_\bq^2(T) := \psi_\bq(\psi_\bq(T))$ and $\psi_\bq^3(T) := \psi_\bq(\psi_\bq^2(T))$.

\begin{prop}
\label{PropFSigmaDec}
Let $T$ be a triangle such that $\sigma_\bq(\psi_\bq^3(T))\geq 5$. Then $\sigma_\bq(\psi_\bq^3(T))\leq 0.69 \sigma_\bq(T)$.
\end{prop}

\proof
The monotonicity of $\sigma_\bq$ established in Proposition \iref{PropSigmaDec} implies that 
$$
5 \leq \sigma_\bq(\psi_\bq^3(T)) \leq \sigma_\bq(\psi_\bq^2(T)) \leq \sigma_\bq(\psi_\bq(T)).
$$
Combining this with inequality \iref{sigma58rho} we obtain 
$$
8\leq \min\{\rho_\bq(\psi_\bq^2(T)), \ \rho_\bq(\psi_\bq(T)), \ \rho_\bq(T)\}.
$$
According to inequality \iref{rhodec}, if a triangle $S$ obeys $\rho_\bq(S)\geq 4$, then 
$
\frac 1 2 \left(1+\frac{16}{\rho_\bq^2(S)}\right)\leq 1
$
and therefore $\rho_\bq(\psi_\bq(S)) \leq \rho_\bq(S)$. We can apply
this to $S = \psi_\bq(T)$ and $S=T$, therefore obtaining
\be
\label{eqMonot}
\rho_\bq(\psi_\bq^2(T)) \leq \rho_\bq(\psi_\bq(T))  \leq \rho_\bq(T).
\ee
We now remark that inequality \iref{rhodec} is equivalent to 
$
(\rho_\bq(S)- \rho_\bq(\psi_\bq(S)))^2 \geq  \rho_\bq(\psi_\bq(S))^2 - 16,
$
hence
\be
\label{decSqrt}
\rho_\bq(S)\geq \rho_\bq(\psi_\bq(S))+ \sqrt{\rho_\bq(\psi_\bq(S))^2-16}
\ee
provided that $\rho_\bq(S) \geq \rho_\bq(\psi_\bq(S))$. 
Applying this to $S=\psi_\bq(T)$ and recalling that $\rho_\bq(\psi_\bq^2(T))\geq 8$ we obtain 
$$
\rho_\bq( \psi_\bq(T))\geq 8+ \sqrt{8^2-16} \geq 14.9.
$$
Applying again \iref{decSqrt} to $S = T$ we obtain 
$$
\rho_\bq(T)\geq 14.9+\sqrt{14.9^2-16} \geq  29.3.
$$
Using \iref{rhodec}, it follows that
$$
\frac{\rho(\psi_\bq^3(T))}{\rho(T)} \leq \frac 1 8 \left(1+\frac{16}{\rho_\bq^2(\psi_\bq^2(T))}\right)\left(1+\frac{16}{\rho_\bq^2(\psi_\bq(T))}\right)\left(1+\frac{16}{\rho_\bq^2(T)}\right)\leq 0.171.
$$
Eventually, the inequalities \iref{equivnondeg} imply that
$$
2 \sigma_\bq(\psi_\bq^3(T))\leq \rho_\bq(\psi_\bq^3(T)) \leq 0.171 \rho_\bq(T) \leq 0.171 (8\sigma_\bq(T))
$$
which concludes the proof.\sq
\nl
An immediate consequence of Propositions \ref{PropSigmaDec} and \ref{PropFSigmaDec} is the following.

\begin{corollary} \label{PropControlSigma2}
If $(T_i)_{i=1}^{8}$ are the eight children obtained from three successive refinement procedures 
from $T$ for the function $q$, then
\begin{itemize}
	\item for all $i$, $\sigma_\bq(T_i)\leq \sigma_\bq(T)$,
	\item there exists $i$ such that $\sigma_\bq(T_i)\leq 0.69\sigma_\bq(T)$ or $\sigma_\bq(T_i)\leq 5$.\\
\end{itemize}
\end{corollary}
We are now ready to prove that most triangles tend to adopt an optimal aspect ratio
as one iterates the refinement procedure.

\begin{theorem} \label{CorolStabilise}
Let $T$ be a triangle, and $\bq$ a positive definite 
quadratic function. Let  $k = \frac{\ln\sigma_\bq(T)-\ln 5}{-\ln(0.69)}$.
Then after $n$ applications of the refinement procedure starting from $T$,
at most $C n^k 7^{n/3}$ of the $2^n$ generated triangles satisfy $\sigma_\bq(S)\geq 5$, 
where $C$ is an absolute constant. Therefore the proportion 
of such triangles tends exponentially fast to $0$ as $n\to +\infty$.
\end{theorem}

\proof
If we prove the proposition for $n$ multiple of $3$, then it will hold for all $n$ (with a larger constant) since $\sigma_\bq$ decreases at each refinement step. We now assume that $n=3m$, and consider the octree with root $T$ obtained by only considering the 
triangles of generation $3 i$ for $i=0,\cdots,n$. 

According to 
Corollary \ref{PropControlSigma2}, for each node of this tree, 
one of its eight children either checks $\sigma_\bq\leq 5$
or has its non-degeneracy measure diminished by a factor $\theta:=0.69$.
We remark that if $\sigma_\bq$ is diminished at least
$k$ times on the path going from the root $T$ to a leaf $S$,
 then $\sigma_\bq(S)\leq 5$.
As a consequence, the number $N(m)$ of triangles $S$
which are such that $\sigma_\bq(S)> 5$ within the generation level $n=3m$ 
is bounded by the number of words in an eight letters 
alphabet $\{a_1,\cdots,a_8\}$ with length $m$ 
and that use the letter $a_8$ at most $k$ times, namely
$$
N(m)\leq \sum_{l=0}^k \binom m l 7^{m-l}
\leq  C m^k 7^m,
$$
which is the announced result.
\sq

The fact that most triangles tend to adopt an
optimal aspect ratio as one iterates the refinement procedure
is a first hint that the approximation
error in the greedy algorithm might satisfy the estimate
\iref{aniser} corresponding to an optimal triangulation.
The following result shows that this is indeed the case,
when this algorithm is applied
on a triangular domain $\Omega$ to a quadratic function $q$
with positive definite associated quadratic form $\bq$. 
The extension of this result to more general $C^2$ 
convex functions on polygonal domains requires a more involved analysis
based on local perturbation arguments
and is the object of the next section. 

\begin{corollary}
\label{corolQuadApprox}
Let $\Omega$ be a triangle, and let $q$ be a quadratic function with positive definite associated quadratic form $\bq$. 
Let $q_N$ be the approximant of $q$ on $\Omega$ obtained by the greedy algorithm 
for the $L^p$ metric, using the $L^1$ decision function \iref{optil1}.
Then 
$$
\limsup_{N\to \infty} N \|q-q_N\|_{L^p(\Omega)}  \leq C \|\sqrt{\det(\bq)}\|_{L^\tau(\Omega)},
$$
where $\frac 1 \tau=\frac 1 p +1$ and where the constant $C$ depends only on 
on the choice of the approximation operator $\cA_T$ used in the definition of the approximant.  
\end{corollary}

\proof
For any triangle $T$, quadratic function $q \in \Pi_2$, and exponent $p$, let 
$$
e'_T(q)_p := \inf_{\pi\in \Pi_1} \|q-\pi\|_{L^p(T)}
$$
be the error of best approximation of $q$ on $T$.
Let $\TO$ be a fixed triangle of area $1$, then for any $q\in \Pi_2$ and $1\leq p\leq \infty$ one has 
$$
e'_\TO(q)_1 \leq e'_\TO(q)_p \leq e_\TO(q)_p \leq e_\TO(q)_\infty.
$$
Furthermore, $e'_\TO(\cdot)_1$ and $e_\TO(\cdot)_\infty$ are semi norms on the finite dimensional space $\Pi_2$ which vanish precisely on the same subspace of $\Pi_2$, namely $\Pi_1$. Hence these semi-norms are equivalent.
It follows that  
\be
\label{eqAOpt}
c_0 \, e_\TO(q)_p \leq e'_\TO(q)_p \leq e_\TO(q)_p
\ee
where $c_0$ is independent $q\in \Pi_2$ and of $p\geq 1$.
Using the invariance property \iref{commuterror}
we find that \iref{eqAOpt} holds for any triangle $T$ in place of $\TO$
with the same constant $c_0$. We also define for any triangulation $\cT$,
$$
e_\cT(f)_p^p := \sum_{T \in \cT} e_T(f)_p^p \quad \text{ and } \quad e'_\cT(f)_p^p := \sum_{T \in \cT} e'_T(f)_p^p,
$$
and we remark that 
$
c_0 \, e_\cT(q)_p \leq e'_\cT(q)_p \leq e_\cT(q)_p.
$
For each $n$, we denote by $\cT_n^u$ the triangulation of $\Omega$ produced by $n$ 
successive refinements based on the $L^1$ decision function \iref{optil1}
for the quadratic function $q$ of interest (note that $\#(\cT_n^u)=2^n$). We also define 
$
\cT_n^\sigma := \{T\in \cT_n^u\sep \sigma_{\bq}(T) > 5\}.
$
Therefore  $\sigma_{\bq}(T) \leq 5$ if $T\notin\cT_n^\sigma$, and
on the other hand we know from Proposition \ref{PropSigmaDec} that 
$\sigma_{\bq}(T) \leq \sigma_{\bq}(\Omega)$ for any $T\in \cT_n^u$.
It follows from Proposition \ref{propequiverrornondeg1} that 
$$
\begin{array}{ll}
e_{\cT_n^u}(q)_p & \leq C_0\(\sum_{T\in \cT_n^u} (\sigma_{\bq}(T)  |T|^{\frac 1 \tau}\sqrt {\det \bq} )^p\)^{\frac 1 p}\\
& \leq C_0 \( 5^p\times  2^n + \sigma_{\bq}(\Omega)^p\#(\cT_n^\sigma)\)^{\frac 1 p} \left(\frac{|\Omega|}{2^n}\right)^{\frac 1 \tau} \sqrt {\det \bq},
\end{array}
$$
where $C_0$ is the constant in \iref{eqEQSigma}. According to Theorem \ref{CorolStabilise},
we know that 
$$
\lim_{n\to +\infty} 2^{-n}\#(\cT_n^\sigma)=0.
$$
Hence 
$$
\limsup_{n \to \infty} 2^n e_{\cT_n^u}(q)_p \leq 5 C_0 \, |\Omega|^{\frac 1 \tau} \sqrt{\det \bq} = 5 C_0 \|\sqrt {\det \bq}\|_{L^\tau(\Omega)}.
$$
We now denote by $\cT^g_n$ the triangulation generated by the greedy procedure with stopping criterion based on the error $\eta_n := C_0^{-1}2^{-\frac n \tau}   \|\sqrt {\det \bq}\|_{L^\tau(\Omega)}$. 
It follows from \iref{eqEQSigma} that for all $T\in \cT_k^u$ with $k\leq n$, one has
$$
e_T(q)_p\geq C_0^{-1}\sigma_{\bq}(T) \|\sqrt {\det \bq}\|_{L^\tau(T)}\geq 
C_0^{-1}2^{-\frac k \tau}\|\sqrt {\det \bq} \|_{L^\tau(\Omega)}\geq \eta_n,
$$
where we used that $|T| = 2^{-k} |\Omega|$ and that the minimal value of $\sigma_\bq$ is $1$.
This shows that $\cT_g^n$ is a refinement of $\cT_n^u$. Furthermore any triangle  $T \in \cT_n^u$ has at most $2^{k(T)}$ children in $\cT_n^g$, where $k(T)$ is the smallest integer such that 
$$
\eta_n \geq  C_0 \, 2^{-\frac{n+k(T)} \tau } \, \sigma_{\bq}(T)  \|\sqrt {\det \bq}\|_{L^\tau(\Omega)}. 
$$ 
Since $\frac 1 2 \leq \tau\leq 1$ we obtain $2^{k(T)} \leq 2^{\frac{k(T)} \tau} \leq 2^{\frac 1 \tau}C_0^2 \sigma_{q}(T)
\leq 4C_0^2 \sigma_{q}(T)$. Hence 
$$
\#(\cT_n^g) \leq 4  C_0^2 \sum_{T \in \cT_n^u} \sigma_{q}(T) \leq 4 C_0^2 \left(5 \times 2^n + \sigma_{\bq}(\Omega) \#(\cT_n^\sigma)\right) = C_1 2^n (1+\ve_n),
$$
where $C_1 = 20\,C_0^2$ and $\ve_n \to 0$ as $n \to \infty$.
If $\cT_N$ is the triangulation generated after $N$ steps
of the greedy algorithm, then there exists $n\geq 0$ such that $\cT_N$ is a refinement of $\cT_n^g$ (hence a refinement of $\cT_n^u$) and $\cT_{n+1}^g$ is a refinement of $\cT_N$. It follows that $\#(\cT_N) \leq \#(\cT_{n+1}^g) \leq C_1 2^{n+1}(1+ \ve_{n+1})$, and 
$$
c_0 \, e_{\cT_N}(q)_p \leq e'_{\cT_N}(q)_p  \leq e'_{\cT_n^u} (q)_p \leq e_{\cT_n^u} (q)_p,
$$
where we have used the fact that $e'_{\cT}(f)_p\leq e'_{\t\cT}(f)_p$ whenever $\cT$ is a refinement
of $\t\cT$. Eventually,
$$
\limsup_{N\to \infty} Ne_{\cT_N}(q)_p\leq \limsup_{n \to \infty} \frac {C_1}{c_0} 2^{n+1}(1+ \ve_{n+1}) e_{\cT_n^u}(q)\leq \frac {10C_0C_1}{c_0} \|\sqrt{\det q}\|_{L^\tau(\Omega)},
$$
which concludes the proof.
\sq

\section{The case of strictly convex functions}

The goal of this section is to prove that the approximation
error in the greedy algorithm applied to a $C^2$ function $f$ satisfies the estimate
\iref{aniser} corresponding to an optimal triangulation.
Our main result is so far limited to 
the case where $f$ is strictly convex.

\begin{theorem}
\label{optitheo}
Let $f\in C^2(\overline \Omega)$ be such that 
$$
d^2f(x)\geq m I,\; \text{ for all } x\in \Omega
$$ 
for some arbitrary but fixed $m >0$ independent of $x$. Let $f_N$
be the approximant obtained by the greedy algorithm 
for the $L^p$ metric, using the $L^1$ decision function \iref{optil1}.
Then 
\be
\limsup_{N\to \infty} N \|f-f_N\|_{L^p}  \leq C \|\sqrt{\det(d^2f)}\|_{L^\tau},
\label{aniser1}
\ee
where $\frac 1 \tau=\frac 1 p +1$ and where $C$
is an absolute constant (i.e. independent of $p$, $f$ and $m$).  
\end{theorem}
Equation \iref{aniser1} can be rephrased as follows : there exists a sequence $\ve_N(f)$ such that $\ve_N(f)\to 0$ as $N \to \infty$ and 
$$
 \|f-f_N\|_{L^p}  \leq \(C \|\sqrt{\det(d^2f)}\|_{L^\tau}+\ve_N(f)\) N^{-1}. 
$$
Note also that since $\|\sqrt{\det(d^2f)}\|_{L^\tau}>0$, 
there exists $N_0(f)$ such that  $\|f-f_N\|_{L^p}  \leq 2C \|\sqrt{\det(d^2f)}\|_{L^\tau} N^{-1}$ for all $N\geq N_0(f)$.
It should be stressed hard that $N_0(f)$ can be arbitrarily large depending on the function $f$. Intuitively,
this means that when $f$ has very large hessian at certain point, it takes more iterations
for the algorithm to generate triangles with a good aspect ratio. 
The extension of this result to strictly concave functions is immediate
by a change of sign. Its extension to arbitrary $C^2$ functions
is so far incomplete, as it is explained in the end of 
the introduction.
The proof of Theorem \ref{optitheo} 
uses the fact that a strictly convex $C^2$ function is
{\it locally} close to a quadratic function with positive definite hessian,
which allows us to exploit the results obtained in \S 3 for 
these particular functions.

\subsection{A perturbation result}
\label{section:Perturb}

We consider a triangle $T$,  a function $f\in C^2(T)$, a convex quadratic function $q$ and $\mu>0$ such that on $T$
\be
\label{ineqfq}
d^2 q\leq d^2 f\leq (1+\mu) \ d^2 q.
\ee
It follows that
$
\det (d^2 q)\leq \det (d^2 f)\leq \det ((1+\mu) d^2 q ) = (1+\mu)^2 \det (d^2 q).
$
Since $\det (d^2 q) = 4 \det(\bq)$, we obtain
\be
\label{ineqdetfq}
2 \|\sqrt{\det \bq}\|_{L^\tau(T)}\leq \nsdf_{L^\tau(T)} \leq 2(1+\mu)  \|\sqrt{\det \bq}\|_{L^\tau(T)}.
\ee
The following Lemma shows how the local errors associated to $f$ and $q$ are close
\begin{prop}
\label{PropErC2}
The exists a constant $C_e>0$, depending only on the operator $\cA_T$ such that 
\be
\label{ineqefq}
(1-C_e\mu) e_T(q)_p \leq e_T(f)_p \leq (1+C_e\mu) e_T(q)_p.
\ee
\end{prop}

\proof
It follows from inequality \iref{ineqfq} that the functions $f-q$ and $(1+\mu)q -f$ are convex, hence
$$
I_T (f-q) - (f-q) \geq 0 \ \text{ and } \ I_T ((1+\mu)q-f) - ((1+\mu)q-f)\geq 0
$$
on the triangle $T$.
We therefore obtain
$$
0 \leq (I_T f-f) - (I_T q-q) \leq \mu (I_T q-q).
$$
There exists a constant $C_0>0$ depending only on $\cA_T$ such that for any $h\in C^0(T)$, 
$$
e_T(h)_p \leq C_0|T|^{\frac 1 p} \|h\|_{L^\infty(T)}.
$$
Furthermore according to Proposition \ref{propequiverrornondeg} there exists a constant $C_1>0$ depending only on $\cA_T$ such that 
$$
|T|^{1/p} \|q-I_T q\|_{L^\infty (T)} \leq C_1  e_T(q)_p.
$$
Hence 
\begin{eqnarray*}
|e_T(f)_p - e_T (q)_p| &\leq & e_T(f-q)_p\\
& = &e_T((I_T f-f) - (I_T q-q))_p \\
 & \leq & C_0 |T|^{\frac 1 p}\|(I_T f-f) - (I_T q-q)\|_{L^\infty (T)}\\
 & \leq & C_0|T|^{\frac 1 p} \|\mu (I_T q-q)\|_{L^\infty (T)}\\
 & \leq & C_0 C_1 \mu e_T(q)_p
\end{eqnarray*}
This concludes the proof of this Lemma, with $C_e=C_0C_1$.
\sq
\nl
Note that using Proposition {\rm \ref{propequiverrornondeg1}},
and assuming that $\mu\leq c_e:=\frac 1{2C_e}$, we have with
$\frac 1 \tau := 1+\frac 1 p$,
\be
\label{eqErrorFQ}
e_T(f)_p \sim e_T(q)_p \sim  \sigma_\bq(T)  \|\sqrt{|\det \bq|}\|_{L^\tau(T)}\sim \sigma_\bq(T) \nsdf_{L^\tau(T)},
\ee
with absolute constants in the equivalence.

We next study the behavior of the decision function $e\mapsto d_T(e,f)$.
For this purpose, we introduce
the following definition.

\begin{definition} \label{def:dg}
Let $T$ be a triangle with edges $a,b,c$. 
A $\delta$-near longest edge bisection with respect to 
the $\bq$-metric  is a 
bisection of any edge $e\in \{a,b,c\}$ such that 
$$
\bq(e)\geq (1-\delta)\max\{\bq(a), \bq(b), \bq(c)\}
$$
\end{definition}

\begin{proposition} 
Assume that $f$ and $q$ satisfy \iref{ineqfq}. Then, the bisection of $T$ prescribed by 
the decision function $e\mapsto d_T(e,f)$ is a $\mu$-near longest edge bisection
for the $\bq$-metric.
\end{proposition}

\proof
It follows directly from Equation \iref{Dint01} that for any edge $e$ of $T$,
$$
D_T(e,q)\leq D_T(e,f) \leq D_T(e,(1+\mu) q),
$$
hence we obtain using \iref{eqDq}
\be
\label{eqDfq}
\frac {|T|} {12} \bq(e) \leq D_T(e,f) \leq (1+\mu) \frac {|T|} {12} \bq(e).
\ee
Therefore the bisection of $T$ prescribed by 
the decision function $e\mapsto d_T(e,f)$ selects an $e$ such that
$$
(1+\mu) \bq(e) \geq \max\{\bq(a), \bq(b), \bq(c)\}.
$$
It is therefore a $\delta$-near longest edge bisection
for the $\bq$-metric with $\delta = \frac \mu {1+\mu}\leq \mu$ and therefore
also a $\mu$-near longest edge bisection. 
\sq

In the rest of this section,
we analyze the difference between a longest edge bisection in the $\bq$-metric and a $\delta$-near longest edge bisection. For that purpose we introduce a distance between triangles :
if $T_1, T_2$ are two triangles with edges $a_1,b_1,c_1$ and $a_2,b_2,c_2$ such that 
\be
\label{orderedEdges}
\bq(a_1)\geq \bq(b_1)\geq \bq(c_1) \ \text{ and } \ \bq(a_2)\geq \bq(b_2)\geq \bq(c_2),
\ee
we define
$$
\Delta_\bq(T_1,T_2) = \max\{|\bq(a_1)-\bq(a_2)|,|\bq(b_1)-\bq(b_2)|,|\bq(c_1)-\bq(c_2)|\}.
$$
Note that $\Delta_\bq$ is a distance up to rigid transformations.

\begin{lemma}
Let $T_1, T_2$ be two triangles, let $(R_1,U_1)$ and $(R_2,U_2)$ be the two pairs of children 
from the longest edge bisection of $T_1$ in the $\bq$-metric, and a $\delta$-near longest edge bisection of $T_2$ in the $\bq$-metric. 
Then, up to a permutation of the pair of triangles $(R_1,U_1)$,  
$$
\max \{\Delta_\bq(R_1,R_2),\Delta_\bq(U_1,U_2)\}\leq \frac 5 4 \Delta_\bq(T_1,T_2)+ \delta \bq(a_2).
$$
where $a_2$ is the longest edge of $T_2$ in the $\bq$-metric.
\end{lemma}

\proof
We assume that the edges of $T_1$ and $T_2$ are named and ordered as in \iref{orderedEdges}.
Up to a permutation,
$R_1$ and $U_1$ have edge vectors 
$b_1,a_1/2,(c_1-b_1)/2$ and $c_1,a_1/2,(b_1-c_1)/2$.
Two situations
might occur for the pair $(R_2,U_2)$:
\begin{itemize}
\item
$\bq(e) < (1-\delta) \bq(a_2)$ 
for $e=b_2$ and $c_2$. In such a case 
the triangle $T_2$ is bisected towards $a_2$, 
so that up to a permutation,
$R_2$ and $U_2$ have edge vectors 
$b_2,a_2/2,(c_2-b_2)/2$ and $c_2,a_2/2,(b_2-c_2)/2$.
Using that $\bq((c-b)/2)=\bq(c)/2+\bq(b)/2-\bq(a)/4$ when $a+b+c=0$, 
it clearly follows that
$$
\max \{\Delta_\bq(R_1,R_2),\Delta_\bq(U_1,U_2)\}\leq \frac 5 4 \Delta_\bq(T_1,T_2).
$$
\item
$\bq(e) \geq (1-\delta) \bq(a_2)$ 
for some $e=b_2$ or $c_2$. In such a case $T_2$ may be bisected say towards $b_2$, 
so that up to a permutation,
$R_2$ and $U_2$ have edge vectors 
$a_2,b_2/2,(c_2-a_2)/2$ and $c_2,b_2/2,(b_2-c_2)/2$.
But since $|\bq(b_2)-\bq(a_2)|\leq \delta \bq(a_2)$,
we obtain that
\be
\label{distPer}
\max \{\Delta_\bq(R_1,R_2),\Delta_\bq(U_1,U_2)\}\leq \frac 5 4 \Delta_\bq(T_1,T_2)+ \delta \bq(a_2).
\ee
\end{itemize}
\sq

We now introduce a perturbed version of the estimates 
describing the decay of the non-degeneracy measure
which were obtained in Proposition \ref{PropSigmaDec} and Corollary \ref{PropControlSigma2}.

\begin{prop}
\label{PropControlSigmaPer}
If $(T_i)_{i=1}^2$ are the two children obtained from a refinement of a triangle $T$ in which a $\delta$-near longest edge bisection in the $\bq$-metric is selected, then
\be
\label{eqSigmaDecPer}
\max\{\sigma_\bq(T_1),\sigma_\bq(T_2)\} \leq (1+4\delta) \sigma_\bq(T).
\ee
If $(T_i)_{i=1}^{8}$ are the eight children of a triangle $T$ 
obtained from three successive refinements in which a $\delta$-near longest edge bisection in the $\bq$-metric is selected, then
\begin{itemize}
       	\item for all $i$, $\sigma_\bq(T_i)\leq \sigma_\bq(T) (1+C_2\delta)$,
	\item there exists $i$ such that $\sigma_\bq(T_i)\leq 0.69 \, \sigma_\bq(T) (1+C_2\delta)$ or $\sigma_\bq(T_i)\leq M$,
	\end{itemize}
where $C_2 = \frac{61} 4$ and $M = 5(1+C_2 \delta)$. 
\end{prop}

\proof
We first prove \iref{eqSigmaDecPer}, and for that purpose we introduce the two children $T'_1, T'_2$  obtained by bisecting the longest edge of $T$ in the $\bq$-metric.
If follows from \iref{distPer} that, up to a permutation of the pair $(T'_1,T'_2)$, 
$$
\max \{\Delta_\bq(T_1,T'_1),\Delta_\bq(T_2,T'_2)\}\leq \delta \bq(a),
$$
where $a$ is the longest edge of $T$ in the $\bq$-metric.
Hence 
\be
\label{eqSigmaDelta}
|\sigma_\bq(T_i)-\sigma_\bq(T'_i)| \leq \frac {2 \Delta_\bq(T_i,T'_i)}{4|T_i|\sqrt{\det(\bq)}}
\leq 2 \delta \frac {\bq(a)}{4|T_i|\sqrt{\det(\bq)}}
\leq 4 \delta \sigma_\bq(T).
\ee
We know from Proposition \ref{PropSigmaDec} that $\max\{\sigma_\bq(T'_1), \sigma_\bq(T'_2)\}\leq \sigma_\bq(T)$. Combining this point with \iref{eqSigmaDelta} we conclude the proof of \iref{eqSigmaDecPer}.\\

We now turn to proof of the second part of the proposition and for that purpose we introduce the eight children $(T'_i)_{i=1}^8$ obtained from three successive refinements of $T$ in which the longest edge in the $\bq$-metric is selected. 
Iterating \iref{distPer}, we find that, up to a permutation of the triangles $(T'_i)_{i=1}^8$, one has 
$$
\max_{i=1,\cdots,8}\Delta_\bq(T_i,T'_i) \leq 
\left(1+\frac 5 4+\left(\frac 5 4\right)^2\right) \delta \bq(a)=\frac {61}{16} \delta \bq(a) = \frac {C_2 \delta} 4 q(a),
$$
where, again, $a$ is the longest edge of $T$ in the $\bq$-metric.
Repeating the argument \iref{eqSigmaDelta} we find that 
\be
\label{eqSigmaDelta2}
\max_{i=1,\cdots,8} |\sigma_\bq(T_i)-\sigma_\bq(T'_i)|\leq C_2\delta \sigma_\bq(T).
\ee
We know from Corollary \ref{PropControlSigma2}
that $\sigma_\bq(T_i')\leq \sigma_\bq(T)$ for all $i$
and that there exists $i$ such that either
$\sigma_\bq(T_i')\leq 0.69 \, \sigma_\bq(T)$ or $\sigma_\bq(T_i')\leq 5$.
Combining this point with \iref{eqSigmaDelta2} we conclude the proof of the proposition.
\sq

\subsection{Local optimality}
\label{section:OptQuad}

Our next step towards the proof of Theorem \ref{optitheo} is
to show that the triangulation produced by the greedy algorithm
is locally optimal in the following sense: 
if the refinement procedure for the function $f$ 
produces a triangle $T\in\cD$ on which $f$ is close
enough to a quadratic function $q$, then the triangles 
which are generated from the refinement of $T$ tend to adopt
an optimal aspect ratio in the $\bq$-metric, 
and a local version of the optimal
estimate \iref{aniser} holds on $T$. 

We first prove that most triangles adopt an optimal aspect ratio
as we iterate the refinement procedure. Our goal is thus to 
obtain a result similar to Theorem \ref{CorolStabilise} 
which was restricted to quadratic functions. However, 
due to the perturbations by $C_2\mu$
that appear in Proposition \ref{PropControlSigmaPer},
the formulation will be slightly different, yet sufficient for
our purposes: we shall prove that the measure of non-degeneracy
becomes bounded by an absolute constant in an average sense,
as we iterate the refinement procedure.

As in the previous section, we
assume that $f$ and $q$ satisfy \iref{ineqfq}. For any $T$, we define $\cT_n^u(T)$ the triangulation of $T$ 
which is built by iteratively applying the refinement procedure for the function $f$
to {\it all generated triangles} up to $3n$ generation levels.
Note that 
$$
\#(\cT_n^u(T))=2^{3n}\;\; {\rm and}\;\; |T'|=2^{-3n}|T|,\;\; T'\in \cT_n^u(T).
$$
For $r>0$, we define
the average $r$-th power of the measure of non-degeneracy of the $2^{3n}$
triangles obtained from $T$ after $3n$ iterations by
$$
\o {\sigma^r_\bq(n)} = \frac 1 {2^{3n}} \sum_{T'\in \cT_n^u(T) } \sigma^r_\bq(T').
$$
We also define
$$
\gamma(r,\mu):= \frac 1 8 \(0.69(1+C_2\mu)\)^r+\frac 7 8 (1+C_2\mu)^r,
$$
where $C_2$ is the constant in Proposition \ref{PropControlSigmaPer}. 
Note that for any $r>0$, the function $\gamma(r,\cdot)$ is continuous and increasing, and that $0<\gamma(r,0)<1$. Hence  for any $r>0$, there
exists $\mu(r)>0$
and $0<\gamma(r)<1$ such that $\gamma(r,\mu)\leq \gamma(r)$, if $0<\mu<\mu(r)$. 

\begin{prop} \label{CorolGamma}
Assume that $f$ and $q$ satisfy \iref{ineqfq} with $0<\mu\leq \mu(r)$. We then have
$$
\o {\sigma^r_\bq(n)} \leq \sigma^r_\bq(T) \gamma(r)^n+\frac {M^r} {8(1-\gamma(r))},
$$
where $M$ is the constant in Proposition \ref{PropControlSigmaPer}.
Therefore 
$$
\o {\sigma^r_\bq(n)}\leq C_3:= 1+\frac {M^r} {8(1-\gamma(r))},
$$ 
if 
$2^{3n}  \geq 8 \sigma_\bq(T)^{\lambda}$ with $\lambda := \frac{3 r \ln 2}{-\ln \gamma(r)}$.
\end{prop}

\proof
Let us use the notations $u= 0.69(1+C_2\mu)$ and $v=(1+C_2\mu)$.
According to Proposition \ref{PropControlSigmaPer}, we
have 
$$
\o {\sigma^r_\bq(n)}\leq  \mE(\sigma_n^r),
$$
where $\mE$ is the expectation operator and $\sigma_n$ is the Markov chain with value in $[1,+\infty[$ defined by
\begin{itemize}
	\item $\sigma_{n+1} = \max\{\sigma_n u, M\}$ with probability $\alpha:=\frac 1 8$,
	\item $\sigma_{n+1} = \sigma_n v$ with probability $\beta:=\frac 7 8$,
	\item $\sigma_0:=\sigma_\bq(T_0)$ with probability $1$.
\end{itemize}
Denoting by $\mu_n$ the probability distribution of $\sigma_n$, we have
\begin{eqnarray*}
\mE(\sigma_{n+1}^r) &=& \int_1^\infty \sigma^r d\mu_{n+1}(\sigma)  \\
&= &\int_1^\infty \left( \alpha (\max\{u\sigma,M\})^r+\beta (v\sigma)^r \right)d\mu_n(\sigma)\\
&=& \alpha M^r \int_1^{M/u} d\mu_n(\sigma) + \alpha u^r \int_{M/u}^\infty \sigma^r d\mu_n(\sigma)
+\beta v^r \int_1^{+\infty}\sigma^r d\mu_n(\sigma)\\
&\leq &\alpha M^r+(\alpha u^r + \beta v^r) \mE(\sigma_n^r)\\
&\leq &
\alpha M^r+\gamma(r) \mE(\sigma_n^r)
\end{eqnarray*}By iteration, it follows that
$$
\mE(\sigma_n^r)\leq \mE(\sigma_0^r)\gamma(r)^n+\frac{ \alpha M^r}{1-\gamma(r)},
$$ 
which gives the result.
\sq

Our next goal is to show that the greedy algorithm initialized from $T$
generates a triangulation which is a refinement 
of $\cT_n^u(T)$ and therefore more accurate, yet with a similar
amount of triangles. To this end, we apply the greedy 
algorithm with root $T$ and stopping criterion given
by the local error 
$$
\eta:=\min_{T'\in \cT_n^u(T)} e_{T'}(f)_p.
$$
Therefore $T'$ is splitted if and only if $e_{T'}(f)_p> \eta$.
We denote by $\cT_{N}(T)$ the resulting triangulation
where $N$ is its cardinality. From the definition of the 
stopping criterion, it is clear that $\cT_{N}(T)$ is a refinement 
of $\cT_n^u(T)$. 

\begin{prop} \label{PropQuadLP}
Assume that $f$ and $q$ satisfy \iref{ineqfq} with $\mu \leq \frac {1}{8}$,
and define
$r_0:=\frac{\ln 2}{\ln 4-\ln 3}>0$.
We then have
$$
N\leq C_42^{3n}\o {\sigma^{r_0}_\bq(n)},
$$
where $C_4$ is an absolute constant.
Assuming in addition that $\mu\leq \mu(r_0)$ as in Proposition \ref{CorolGamma},
we obtain that
$$
N\leq C_52^{3n},
$$
if $2^{3n} \geq 8 \sigma_\bq(T)^{\lambda}$ with $\lambda := \frac{3 r_0\ln 2}{-\ln \gamma(r_0)}$,
and where $C_5=C_3C_4$.
\end{prop}

\proof
Let $T_1$ be a triangle in $\cT_{n}^u(T)$ and $T_2$ a triangle in $\cT_{N}(T)$
such that $T_2\subset T_1$. We shall give a bound on the number of splits
$k$ which were applied between $T_1$ and $T_2$, i.e. such that
$|T_2|=2^{-k}|T_1|$. We first remark that according to Proposition 
\ref{propequiverrornondeg1} and \iref{eqErrorFQ}, we have
$$
\eta\geq c\min_{T'\in \cT_n^u(T)}|T'|^{1+\frac 1 p} \sigma_\bq(T')  \sqrt{\det \bq}
\geq c|T_1|^{1+\frac 1 p} \sqrt{\det \bq},
$$
where $c$ is an absolute constant.
On the other hand, using both
Proposition \ref{PropErC2} and Proposition \ref{PropControlSigmaPer}, we obtain
$$
\begin{array}{ll}
e_{T_2}(f)_\bq 
&\leq C|T_2|^{1+\frac 1 p} \sigma_\bq(T_2)  \sqrt{\det \bq} \\
&=C|T_1|^{1+\frac 1 p} 2^{-k(1+\frac 1 p)}\sigma_\bq(T_2)  \sqrt{\det \bq} \\
&\leq C|T_1|^{1+\frac 1 p} \sigma_\bq(T_1)\(2^{-(1+\frac 1 p)}(1 + 4\mu)\)^{k} \sqrt{\det \bq}.\\
&\leq \frac C c \sigma_\bq(T_1)\(\frac {1+4\mu} 2\)^{k} \eta \\
&\leq \frac C c \sigma_\bq(T_1)(\frac 3 4)^{k} \eta,
\end{array}
$$
where $C$ is an absolute constant.
Therefore we see that $k$ is at most the smallest integer 
such that $\frac C c \sigma_\bq(T_1)(\frac 3 4)^{k} \leq 1$.
It follows that the total number $n(T_1)$ of triangles $T_2\in \cT_{N}(T)$
which are contained in $T_1$
is bounded by
$$
n(T_1)\leq 2^k \leq 2 \left(\frac C c \sigma_\bq(T_1)\right)^{r_0},
$$
and therefore
$$
N=\sum_{T_1\in \cT^u_n(T)} n(T_1) \leq 2 \left(\frac C c\right)^{r_0} \sum_{T_1\in \cT^u_n(T)} 
\sigma_\bq(T_1)^{r_0} = C_42^{3n}\o {\sigma^{r_0}_\bq(n)},
$$
with $C_4= 2 \left(\frac C c\right)^{r_0}$. The fact that 
$N\leq C_52^{3n}$ when $2^{3n}  \geq 8 \sigma_\bq(T)^{\lambda}$ 
with $\lambda := \frac{3 r_0 \ln 2}{-\ln \gamma(r_0)}$ is an immediate
consequence of Proposition {\rm \ref{CorolGamma}}.
\sq

\subsection{Optimal convergence estimates}

Our last step towards the proof of Theorem
\ref{optitheo} consists in deriving local error estimates for the greedy algorithm.
For $\eta>0$, we 
denote by $f_{\eta}$ the approximant to $f$
obtained by the greedy algorithm with
stopping criterion given
by the local error $\eta$ : a triangle
$T$ is splitted if and only if $e_{T}(f)_p > \eta$.
The resulting triangulation is denoted by
$$
\cT_\eta=\cT_N, \text{ with } N=N(\eta)=\#(\cT_\eta).
$$
For this $N$, we thus have $f_\eta=f_N$.
For a given $T$ generated by the refinement procedure and 
such that $\eta \leq e_T(f)_p$, we also define
$$
\cT_{\eta}(T)=\{T'\subset T\; ; \; T' \in \cT_\eta\}
$$
the triangles in $\cT_\eta$ which are contained in $T$
and
$$
N(T,\eta)=\#(\cT_\eta(T)).
$$
Our next result provides with estimates of the local
error $\|f-f_\eta\|_{L^p(T)}$ and of $N(T,\eta)$ in terms of 
$\eta$, provided that $\mu$ is small enough.

\begin{theorem} 
\label{localerror}
Assume that $f$ and $q$ satisfy \iref{ineqfq} with $\mu\leq c_2:= \min\{ \frac 1 8,\mu(r_0)\}$,
and that $\eta \leq \eta_0$, where
$$
\eta_0=\eta_0(T):=\left(\frac {|T|}{\sigma_\bq(T)^{\lambda}}\right)^{\frac 1 \tau}\sqrt{\det \bq},
$$
with $\lambda := \frac{3 r_0\ln 2}{-\ln \gamma(r_0)}$,
and $\frac 1 \tau=\frac 1 p+1$. Then
\be
\|f-f_\eta\|_{L^p(T)}\leq \eta N(T,\eta)^{\frac 1 p},
\label{localerreta}
\ee
and 
\be
N(T,\eta) \leq C_6\eta^{-\tau} \|\sqrt{\det(d^2f)}\|_{L^\tau(T)}^\tau,
\label{estimNeta}
\ee
where $C_6$ is an absolute constant.
\end{theorem}

\proof
The first estimate is trivial since
$$
\|f-f_\eta\|_{L^p(T)}=\(\sum_{T'\in \cT_{\eta}(T)} e_{T'}(f)_p^p\)^{\frac 1 p}
\leq \(\sum_{T'\in \cT_{\eta}(T)} \eta^p\)^{\frac 1 p} =\eta N(T,\eta)^{\frac 1 p}.
$$
In the case $p=\infty$, we trivially have
$$
\|f-f_\eta\|_{L^\infty(T)}\leq \eta.
$$
For the second estimate, we define $n_0=n_0(T)$ the smallest
positive integer such that 
$2^{3n_0(T)}\geq  8\sigma_\bq(T)^{\lambda}$ with 
$\lambda := \frac{3 r_0\ln 2}{-\ln \gamma(r_0)}$.
For any fixed $n\geq n_0$, we define
$$
\eta_n:=\min_{T'\in \cT_n^u(T)} e_{T'}(f)_p.
$$
We know from Proposition \ref{PropQuadLP} that with the choice $\eta=\eta_n$
\be
N(T,\eta_n)\leq C_5 2^{3n}.
\label{netan}
\ee
On the other hand, we know from Proposition \ref{CorolGamma}, 
that $\o{\sigma_\bq^{r_0}(n)}\leq C_3$,
from which it follows that 
$$
\min_{T'\in\cT_n^u(T)} \sigma_\bq(T') \leq C_3^{\frac 1 {r_0}}.
$$
According to Proposition 
\ref{PropErC2}, we also have 
$$
\eta_n
\leq C\min_{T'\in \cT_n^u(T)}|T'|^{1+\frac 1 p} \sigma_\bq(T')  \sqrt{\det \bq}
\leq C_3^{\frac 1 {r_0}} C  \(\frac {|T|}{2^{3n}}\)^{\frac 1 \tau}\sqrt{\det \bq},
$$
where $C$ is an absolute constant, which also reads
$$
2^{3n}\leq C_3^{\frac \tau {r_0}}C^\tau  \eta_n^{-\tau} |T|\sqrt{\det \bq}^{\,\tau}.
$$
Combining this with \iref{netan}, we have obtained the estimate
$$
N(T,\eta_n)\leq C_5C_3^{\frac \tau {r_0}}C^\tau  \eta_n^{-\tau} |T|\sqrt{\det \bq}^{\,\tau},
$$
which by Proposition \ref{PropErC2} is equivalent to \iref{estimNeta} 
with $\eta=\eta_n$. In order to obtain \iref{estimNeta} for all arbitrary
values of $\eta$, we write that 
$\eta_{n+1}< \eta \leq \eta_n$ for some $n\geq n_0$, then
$$
\begin{array}{ll}
N(T,\eta) &\leq N(T,\eta_{n+1}) \\
& \leq C_5 2^{3(n+1)}\\
& \leq  8C_5C_3^{\frac \tau {r_0}}C^\tau  \eta_n^{-\tau} |T|\sqrt{\det \bq}^{\,\tau}\\
&\leq 8C_5C_3^{\frac \tau {r_0}}C^\tau  \eta^{-\tau} |T|\sqrt{\det \bq}^{\,\tau},
\end{array}
$$
which by Proposition \ref{PropErC2} is equivalent to \iref{estimNeta}.
In the case where $\eta \geq \eta_{n_0}$, we simply write
$$
\begin{array}{ll}
N(T,\eta) & \leq  N(T,\eta_{n_0}) \\
& \leq C_5 2^{3n_0}\\
& \leq 64C_5\sigma_\bq(T)^{\lambda} \\
&=64C_5\eta_0^{-\tau} |T|\sqrt{\det \bq}^{\,\tau} \\
& \leq 64C_5\eta^{-\tau} |T|\sqrt{\det \bq}^{\,\tau},
\end{array}
$$
and we conclude in the same way.
\sq
\nl
We remark that combining the estimates
\iref{localerreta} and \iref{estimNeta} in the above theorem
yields the optimal local convergence estimate
$$
\|f-f_\eta\|_{L^p(T)}\leq C_6^{\frac 1 \tau} \|\sqrt{\det(d^2f)}\|_{L^\tau(T)}N(T,\eta)^{-1}.
$$
In order to obtain the global estimate of Theorem \ref{optitheo}, we 
need to be ensured that after sufficiently many steps of the greedy algorithm, 
the target $f$ can be well approximated by quadratic function $q=q(T)$ 
on each triangle $T$, so that our local results will apply on such triangles.
This is ensured due to the following key result.

\begin{prop} \label{PropDiamZero}
Let $f$ be a $C^2$ function such that 
$d^2 f(x)\geq m I $ 
for some arbitrary but fixed $m >0$ independent of $x$. 
Let $\cT_N$ be the triangulation generated by the greedy algorithm
applied to $f$ using the $L^1$ decision function given by \iref{optil1}.
Then
$$
\lim_{N\to +\infty} \max_{T\in\cT_N}\diam(T)=0,
$$
i.e. the diameter of all triangles tends to $0$.
\end{prop}
\proof
Let $T$ be a triangle with an angle $\theta$ at a vertex $z_0$. The other vertices of $T$ 
can be written as $z_1=z_0+\alpha u$ and $z_2=z_0+ \beta v$ where $\alpha, \beta\in \R_+$ and $u,v\in \R^2$ are unitary. We assume that $\alpha u$ is the longest edge of $T$, hence $\theta\leq \cPi/2$.
Observe that 
$$
\rho(T) := \frac {h_T^2} {|T|}=\frac{\alpha^2} {\frac 1 2\alpha \beta \sin \theta} = \frac {2 \alpha}{\beta \sin \theta},
$$
and
$$ 
|u-v| = 2 \sin\left( \frac \theta 2\right)=\frac {\sin \theta}{\cos(\frac \theta 2)}
=\frac {2\alpha}{\beta\rho(T)\cos(\frac \theta 2)}.
$$
Since $\frac {\sqrt 2} 2 \leq \cos\left(\frac \theta 2\right) $ we thus obtain 
$$
|u-v| \leq  \frac {2\sqrt 2\alpha}{\beta \rho(T)} \leq \frac {3\alpha}{\beta \rho(T)} .
$$
We now set $M := \|d^2 f\|_{L^\infty(\Omega)}$ and for all $\delta>0$ let 
$$
\omega(\delta) := \sup_{z,z'\in \Omega, \|z-z'\|\leq \delta} \|d^2f(z) - d^2f(z')\|.
$$
For $t\in\RR$, we define 
$$
H_t^u:= d^2 f_{z_0 + t u} \ \text{ and } \ H_t^v:= d^2 f_{z_0 + t v}.
$$
and notice that 
$
\|H_t^u - H_t^v\| \leq \omega(t|u-v|).
$
Hence, if $0\leq t\leq \beta$, we have
$$
\|H_t^u-H_t^v\|\leq \omega \left(\frac {3\alpha } {\rho(T)}\right).
$$ 
Furthermore, for all $t$ we have
$$
\begin{array}{ll}
|\<H_t^u \ u, u\> - \<H_t^u \ v, v\>| & 
=|\<H_t^u \ u, u\> - \<H_t^u \ u-(u-v), u-(u-v)\>| \\
& =|2\<H_t^u \ u, u-v\> -\<H_t^u \ (u-v), u-v\>|\\
&\leq 2M |u| |u-v|+ M |u-v|^2 \\
& \leq \frac M {\beta^2} \left( 2\frac {3\alpha \beta}{\rho(T)} + \left(\frac {3\alpha}{\rho(T)}\right)^2 \right).
\end{array}
$$
Applying the identity \iref{Dint01} to the edges $e=\alpha u$ and $\beta v$,
and using a change of variable, we can write
$$
D_T(\alpha u ,f)=\int_\R \min\{t,\alpha-t\}_+ \<H_t^u\, u,u\> dt
\;\; {\rm and}\;\; D_T(\beta v ,f)=\int_\R \min\{t,\beta-t\}_+ \<H_t^v\, v,v\> dt
$$
where we have used the notation $r_+ := \max\{r,0\}$.
Hence, noticing that 
$$
\int_\R \min\{t,\lambda-t\}_+ dt=\int_0^\lambda \min\{t,\lambda-t\}dt = \frac{(\lambda_+)^2} 4,
$$
and using the previous estimates
we obtain
\begin{eqnarray*}
D_T(\alpha u ,f) - D_T(\beta v,f) &=& \int_\R \left(\min\{t,\alpha-t\}_+ \<H_t^u\, u,u\> - \min\{t, \beta-t\}_+ \<H_t^v \, v,v\>\right) dt\\
&= &  \int_\R (\min\{t,\alpha-t\}_+ -  \min\{t, \beta-t\}_+) \<H_t^u \, u,u\> dt\\
& & -\int_\R  \min\{t, \beta-t\}_+ (\<H_t^v \, v,v\>- \<H_t^u\, u,u\>)dt \\
& \geq &m\int_\R (\min\{t,\alpha-t\}_+ -  \min\{t, \beta-t\}_+)dt \\
& & -\int_0^\beta \min\{t, \beta-t\} (|\<H_t^u \ u, u\> - \<H_t^u \ v, v\>| +| \<(H_t^u-H_t^v)v,v\>|)dt \\
&\geq &  m \frac {\alpha^2 -\beta^2} 4  -
\frac M 4 \left( 2\frac {3\alpha\beta}{\rho(T)} + \left(\frac {3\alpha}{\rho(T)}\right)^2 \right)
- \frac{\beta^2} 4 \omega\left(\frac { 3\alpha}{ \rho(T)}\right)\\
& \geq & m \frac {\alpha^2 -\beta^2} 4 - \frac {\alpha^2} 4 \left(
\frac 6 {\rho(T)}+\frac 9 {\rho(T)^2}+ \omega\left(\frac { 3\alpha}{ \rho(T)}\right)\right),
\end{eqnarray*}
where we have used the fact that $\alpha>\beta$ in the last line.
We can therefore write
\be
D_T(\alpha u ,f) - D_T(\beta v,f) \geq \frac m 4 \(\alpha^2\(1-\t \omega\(\frac 1{\rho(T)}\)\)-\beta^2\),
\label{eqDT}
\ee
where we have set
$$ 
\tilde \omega(\delta) := \frac 1 m \( 6 \delta+9 \delta^2+\omega( 3\diam(\Omega) \delta)\).
$$ 
The inequality \iref{eqDT} shows that $d_T(\cdot,f)$ prescribes a
$\t \omega\(\frac 1 {\rho(T)}\)$-near longest edge bisection in the euclidean metric for any triangle $T$. 
Indeed if the smaller edge $\beta v$ was selected, we would necessarily have
$$
|\beta v|^2=\beta^2 \geq \(1-\t \omega\(\frac 1 {\rho(T)}\)\) \alpha^2 = \(1-\t \omega\(\frac 1 {\rho(T)}\)\)
|\alpha u|^2.
$$
Notice that $\tilde \omega(\delta)\to 0$ as $\delta \to 0$.

Since $f$ is strictly convex, there does not exists any triangle $T\subset \Omega$ such that $e_T(f)_p = 0$. 
Let us assume for contradiction that the diameter of the triangles generated by the greedy
algorithm does not tend to zero. Then there 
exists a sequence $(T_i)_{i\geq 0}$ of triangles such that 
$T_{i+1}$ is one of the children of $T_i$, and $h_{T_i}\to d>0$ 
as $i \to \infty$, where $h_T$ denotes the diameter of a triangle $T$. Since $|T_i|\to 0$, this also implies that $\rho(T_i)\to +\infty$
as $i\to \infty$. We can therefore choose $i$ large enough such that
$h_{T_i}^2<\frac 4 3 d^2$ and 
$C_2\tilde \omega\(\frac 1 {\rho(T_j)}\)\leq \frac 1 2$ for all $j\geq i$,
where $C_2$ is the constant in Proposition \ref{PropControlSigmaPer}.
According to this Proposition, we have
$$
\sigma(T_{i+3}) \leq \frac 3 2 \sigma(T_i),
$$
where $\sigma$ stands for $\sigma_{\bq}$ in the euclidean case $\bq=x^2+y^2$.
On the other hand, we have for any triangle $T$,
$
\frac {h_T^2}{8|T|} \leq \sigma(T)\leq \frac {h_T^2}{2|T|},
$
from which it follows that 
$$
h_{T_{i+3}}^2 \leq 4 \frac {|T_{i+3}|\sigma(T_{i+3})}  {|T_{i}|\sigma(T_{i})}
h_{T_i}^2 \leq \frac 3 4
h_{T_i}^2.
$$
Therefore, $h_{T_{i+3}}<d$ 
which is a contradiction. This concludes the proof of
Proposition \ref{PropDiamZero}.
\sq

\noindent
{\bf Proof of Theorem \ref{optitheo}}
Since $f\in C^2$, an immediate consequence of Proposition
\ref{PropDiamZero} is that for
all $\mu>0$, there exists 
$$
N_1:=N_1(f,\mu),
$$
such that for all $T\in \cT_{N_1}$, there exists a quadratic function $q_T$
such that  
$$
d^2 q_T \leq d^2 f\leq (1+\mu) \, d^2 q_T.
$$
Therefore our local results apply on all $T\in \cT_{N_1}$. Specifically,
we choose
$$
N_1:=N_1(f,c_2),
$$ 
with $c_2$ the constant in Theorem \ref{localerror}. We then take 
$$
\eta \leq \eta_0:=\min_{T\in\cT_{N_1}}\left\{e_T(f)_p,\(\frac {|T|}{ \sigma_{\bq_T}(T)^{\lambda}}\)^{\frac 1 \tau}\sqrt{\det \bq_T}\right\}.
$$
We use the notations
$$
f_{\eta}=f_N,\;\; \cT_\eta=\cT_N,\;\, N=N(\eta)=\#(\cT_\eta)=\#(\cT_N),
$$
for the approximants and triangulation
obtained by the greedy algorithm with
stopping criterion given
by the local error $\eta$. Note that $\cT_{\eta}$ is a refinement of $\cT_{N_1}$,
since $\eta\leq \min_{T\in\cT_{N_1}}e_T(f)_p$, and therefore $N\geq N_1$.
We obviously have
$$
\|f-f_N\|_{L^p}\leq \eta N^{\frac 1 p}.
$$
Using Theorem \ref{localerror}, we also have
$$
N=\sum_{T\in \cT_{N_1}}N(T,\eta)\leq C_6\eta^{-\tau} \|\sqrt{\det(d^2f)}\|_{L^\tau(\Omega)}^\tau,
$$
and therefore 
$$
\|f-f_N\| \leq  C_6^{\frac 1 \tau} \|\sqrt{\det(d^2f)}\|_{L^\tau(\Omega)}N^{-1},
$$
which is the claimed estimate. Since we have assumed 
$\eta\leq \eta_0$, this estimate holds for 
$$
N> N_0,
$$
where $N_0$ is largest value of $N$ such that $e_T(f)_p\geq \eta_0$
for at least one $T\in\cT_N$.
\sq

\begin{remark}
In \cite{CDHM} a modification of the algorithm is proposed so that 
its convergence in the $L^p$ norm is ensured
for {\it any} function $f\in L^p(\Omega)$ (or $f\in C(\Omega)$ when $p=\infty$). 
However this modification is not needed 
in the proof of Theorem \ref{optitheo}, due to the assumption
that $f$ is convex. 
\end{remark}

$\;$
\nl
Albert Cohen
\nl
UPMC Univ Paris 06, UMR 7598, Laboratoire Jacques-Louis Lions, F-75005, Paris, France
\nl
CNRS, UMR 7598, Laboratoire Jacques-Louis Lions, F-75005, Paris, France
\nl
cohen@ann.jussieu.fr
\nl
\nl
Jean-Marie Mirebeau
\nl
UPMC Univ Paris 06, UMR 7598, Laboratoire Jacques-Louis Lions, F-75005, Paris, France
\nl
CNRS, UMR 7598, Laboratoire Jacques-Louis Lions, F-75005, Paris, France
\nl
mirebeau@ann.jussieu.fr


\begin {thebibliography} {99}


\bibitem{BBLS} V. Babenko, Y. Babenko, A. Ligun and A. Shumeiko,
{\it On Asymptotical
Behavior of the Optimal Linear Spline Interpolation Error of $C^2$
Functions}, East J. Approx. 12(1), 71--101, 2006.


\bibitem{BDDP} P. Binev, W. Dahmen, R. DeVore and P. Petrushev,
{\it Approximation Classes for Adaptive Methods}, 
Serdica Math. J. 28, 391--416, 2002.

\bibitem{BFGLS} H. Borouchaki, P.J. Frey,  P.L. George, P. Laug and E. Saltel,
{\it Mesh generation and mesh adaptivity: theory, techniques},  in Encyclopedia of computational mechanics, E. Stein, R. de Borst and T.J.R. Hughes ed., John Wiley \& Sons Ltd., 2004.


\bibitem{Cao} W. Cao, {\it On the error of linear interpolation and the orientation, aspect ratio, and
internal angles of a triangle}, SIAM J. Numer. Anal. 43(1), 19--40, 2005.

\bibitem{Chen1} L. Chen, {\it On minimizing the linear interpolation error
of convex quadratic functions}, East Journal of Approximation 14(3), 271--284, 2008.

\bibitem{Chen2} L. Chen, {\it Mesh smoothing schemes based on optimal Delaunay
triangulations}, in {\it 13th International Meshing Roundtable}, 109--120, Williamsburg VA, Sandia National
Laboratories, 2004.

\bibitem{CSX} L. Chen, P. Sun and J. Xu, {\it Optimal anisotropic meshes for
minimizing interpolation error in $L^p$-norm}, Math. of Comp. 76, 179--204, 2007.


\bibitem{CDHM} A. Cohen, N. Dyn, F. Hecht and J.-M. Mirebeau,
{\it Adaptive multiresolution analysis based on anisotropic triangulations}, accepted in Mathematics of Computation 2010.

\bibitem{De} R. DeVore, {\it Nonlinear approximation},
Acta Numerica 51-150, 1998

\bibitem{Mi} J.-M. Mirebeau, {\it Optimally adapted finite elements meshes}, Constructive Approximation, Vol 32 $\text{n}^\text{o}$2, pages 339-383, 2010.

\bibitem{thesisJM} J.-M. Mirebeau, {\it Adaptive and anisotropic finite element approximation : Theory and algorithms}, PhD Thesis, \url{tel.archives-ouvertes.fr/tel-00544243/en/}


\bibitem{Ri} M.C. Rivara, {\it New longest-edge algorithms for the refinement and/or improvement of unstructured triangulations}, 
Int. J. Num. Methods 40, 3313--3324, 1997.


\end{thebibliography}
\end{document}